\documentclass[reqno]{amsart}
\usepackage[usenames,dvipsnames]{color}
\usepackage{amsthm}
\usepackage{amsmath}
\usepackage{amssymb}
\usepackage{amscd}
\usepackage{graphics}
\usepackage{latexsym}
\usepackage{stmaryrd}
\usepackage{empheq}
\usepackage{xcolor}
\usepackage[colorlinks=true,linkcolor=black,
 citecolor=orange,urlcolor=black]{hyperref}
\usepackage{enumitem}
\usepackage[style=trad-abbrv,maxnames=99,maxalphanames=9, isbn=false, giveninits=true, doi=false, url=true]{biblatex}
\bibliography{ref.bib}
\renewbibmacro{in:}{}
\theoremstyle{remark}

\definecolor{titlecol}{named}{BrickRed}
\definecolor{headcol}{named}{Violet}
\definecolor{seccol}{named}{Red}
\definecolor{sseccol}{named}{Bittersweet}
\definecolor{pbcol}{named}{Black}
\definecolor{sncol}{named}{Brown}
\definecolor{acol1}{named}{Red}
\definecolor{acol2}{named}{Apricot}

\def\R{\mathbb{R}}
\def\C{\mathbb{C}}
\def\N{\mathbb{N}}
\def\Z{\mathbb{Z}}

\def\tr{\text{tr}}
\newcommand{\ddc}{dd^c}
\newcommand{\Psef}{\operatorname{Psef}}
\newcommand{\MN}{\operatorname{MN}}
\newcommand{\MK}{\operatorname{MK}}
\newcommand{\NS}{\operatorname{NS}}
\newcommand{\JNull}{\operatorname{Null}_{J}}
\newcommand{\JEnK}{E_{nJ}}
\newcommand{\wto}{\rightharpoonup}
\newcommand{\Psh}{\mathcal{PSH}}
\newcommand{\Eff}{\text{Eff}}
\newcommand{\cab}{c_{\alpha,\beta}}
\newcommand{\reg}{\text{reg}}
\newcommand{\sing}{\text{Sing}}
\newcommand{\supp}{\text{supp}}
\newcommand{\Null}{\operatorname{Null}}

\def\cG{{\mathcal G}}

\def\cI{{\mathcal I}}
\def\cJ{{\mathcal J}}
\def\cK{{\mathcal K}}

\def\cO{{\mathcal O}}

\theoremstyle{plain}
\newtheorem{thm}{Theorem}

\newtheorem{cor}[thm]{Corollary}
\newtheorem{prop}[thm]{Proposition}
\newtheorem{lemma}[thm]{Lemma}
\newtheorem{rmk}[thm]{Remark}

\newtheorem{defn}[thm]{Definition}

\newtheorem{remark}[thm]{Remark}

\title{Semistable $J$-equation on K\"ahler threefold}
\author{Junbang Liu}
\address{Department of Mathematics, The Hong Kong University of Science and Technology, Clear Water Bay, Kowloon, Hong Kong}
\email{junbangliu@ust.hk}

\begin{document}
\maketitle
\begin{abstract}
    For a $3$-dimensional compact K\"ahler manifold $X$ with a pair of K\"ahler class $(\alpha,\beta)$ which is $J$-semistable, we show that the $J$-null locus of $(\alpha,\beta)$ forms a subvariety of $X$. We also give an analyic characterization of the $J$-null locus using K\"ahler currents with analytic singularities analogous to the theorem of Collins-Tosatti \cite{CT15}. As an application, we show the smooth convergence of the $J$-flow off the $J$-null locus. It gives higer order partial regularity for weak solutions of semistable $J$-equations obtained in \cite{M26b}. 
\end{abstract}
\tableofcontents
\section{Introduction}
Let $X$ be a compact K\"ahler manifold of complex dimension $n$. Fix K\"ahler classes $\alpha,\beta$, and choose K\"ahler representatives $\omega\in \alpha, \chi\in \beta$, and set $c_{\alpha,\beta}=n\beta\cdot\alpha^{n-1}/{\alpha^n}$. 
The $J$-equation asks for a K\"ahler metric $\omega_\varphi:=\omega+\ddc \varphi\in \alpha$, satisfying $\tr_{\omega_\varphi}\chi=\cab$, or equivalently, \begin{equation}
\label{eq-J-equation}
n\omega_\varphi^{n-1}\wedge\chi=\cab\omega_\varphi^n.
\end{equation}
It arose in Donaldson's moment-map picture \cite{D99} and in Chen's study of the Mabuchi energy and the $J$-functional \cite{C00}. It's the Euler-Lagrange equation for the $J$-functional, and it's solvability is closely related to lower bounds and coercivity properties of the Mabuchi energy. 
The natural strict cone, or subsolution, condition associated with the $J$-equation is the existence of a K\"ahler form $\widehat\omega\in\alpha$ such that
\begin{equation}\label{eq:intro-strict-cone}
 \cab\widehat\omega^{\,n-1}
 -(n-1)\chi\wedge\widehat\omega^{\,n-2}>0.
\end{equation}
Song-Weinkove proved, using the $J$-flow, that solvability of the $J$-equation is equivalent to the strict cone condition \cite{SW08}.  
Restricting the inequality \eqref{eq:intro-strict-cone} to subvarieties leads to the numerical conditions \begin{equation}\label{eq:intro-strict-numerical}
 \int_V\bigl(\cab\alpha^p-p\beta\alpha^{p-1}\bigr)>0
 \qquad
 (1\leq p=\dim V\leq n-1).
\end{equation}
Lejmi--Sz\'ekelyhidi formulated these inequalities as the expected
numerical criterion for the $J$-equation \cite{LS15}. Collins-Sz\'ekelyhidi proved the equivalence between the solvability and the coercivity of the $J$-functional \cite{CS17}. The
decisive breakthrough was Chen's theorem
\cite[Theorem~1.1]{C21}.  He proved the equivalence between the
solvability and a quantitative uniform
version of the intersection inequalities \eqref{eq:intro-strict-numerical}. For further developments see \cite{DP21,S20}. Fang-Ma later introduced a fully nonlinear equation with differential-form data which simultaneously contains the Monge-Amp\'ere equations, the $J$-equations, and the inverse $\sigma_k$-equations, and obtained the corresponding analytic and numerical solvability criteria \cite{FM24}.

The boundary and unstable regimes have recently received considerable attention. For the K\"ahler surface case, the equation can be formulated as a complex Monge-Amp\`ere equation, see \cite{FLSW14,DMS26,M26a}. Especially, Datar--Mete--Song introduced a minimal-slope and bubbling program for
complex Hessian equations; for the $J$-equation on K\"ahler surfaces
they constructed canonical singular solutions in the unstable case
and studied flow limits and bubbling in symmetric examples
\cite{DMS26}. 
 Khalid and Sj\"ostr\"om Dyrefelt developed
wall--chamber decompositions for generalized Monge--Amp\`ere equations:
under their positivity hypotheses, only finitely many subvarieties
violate the numerical criterion, and these subvarieties satisfy a
rigidity property \cite{KSD26}.  Most recently, Fu proved
the Datar--Mete--Song minimal-slope conjecture, characterizing
$J$-slope semistability by the minimal $J$-slope
\cite{F26}.  Fu also studied the $J$-null locus of a
semistable pair and proved its analyticity for compact K\"ahler
surfaces and compact toric K\"ahler manifolds. These developments make
the geometry of the degeneracy set, rather than merely weak convergence
or smooth solvability, a natural next question.

The purpose of the present paper is to study this question on an
arbitrary compact K\"ahler threefold.  We assume only that the
expressions in \eqref{eq:intro-strict-numerical} are nonnegative and
ask where strict positivity must fail analytically. There is a compelling model for this question in the theory of
positive $(1,1)$-classes.  If $\gamma$ is nef and big, its numerical
null locus is
\[
 \operatorname{Null}(\gamma)
 :=\bigcup_{\substack{V\subset X\ {\rm irreducible}\\
                       \dim V>0,\ 
                       \int_V\gamma^{\dim V}=0}}V.
\]
Collins--Tosatti proved that this set is the non-K\"ahler locus
$E_{\rm nK}(\gamma)$ \cite[Theorem~1.1]{CT15}.  Thus a
purely numerical degeneracy set is exactly the unavoidable singular
set of K\"ahler currents in the class.  This suggests defining the
$J$-null locus by
\[
 \JNull(\alpha,\beta)
 :=\bigcup_{\substack{V\subsetneq X\ {\rm irreducible}\\
                       \dim V=p>0,\\
  \int_V(c\alpha^p-p\beta\alpha^{p-1})=0}}V
\]
and looking for a corresponding analytic locus defined by currents
satisfying the weak $J$-cone condition. We fix some notations now. For an irreducible proper subvariety $V$ of dimension $p$, set $J_{\cab}(V,\alpha,\beta):=\int_V (\cab\alpha^p-p\beta\alpha^{p-1})$. 
The pair $(\alpha,\beta)$ is called \emph{ J-nef} if $ J_{c_{\alpha,\beta}}(V,\alpha,\beta)\geq 0$ for every proper positive-dimension irreducible analytic subvariety $V\subset X$. For the weak sense of the cone condition \eqref{eq:intro-strict-cone},
we adopt its definition from Chen\cite[Definition 3.3]{C21}. If $A$ is a smooth positive $(1,1)$-form, define
\[
 P_\chi(A)(x)
 :=\max_{\substack{H\subset T_x^{1,0}X\\ \dim_{\mathbb C}H=n-1}}
        \tr_{A|H}(\chi|H).
\]
Thus $P_\chi(A)<c$ is equivalent to $cA^{n-1}-(n-1)A^{n-2}\wedge \chi>0$. If $T$ is a positive current, we write $P_\chi(T)\leq c$ if on any open subset $O$ of any coordinate chart, for any K\"ahler form $\chi_0\leq \chi$ with constant coefficients, $T=\ddc\varphi$ on $O$, one has $P_{\chi_0}(\ddc\varphi_\delta)\leq c,$ for any $\delta>0$.  Here $\varphi_\delta$ is the local convolution smoothing of $\varphi$: $\varphi
   _\delta(x)=\int_{\C^n}\varphi(x-y)\delta^{-2n}\rho(|y|/\delta)dvol_y$ for a smooth nonnegative kernel $\rho$. Now we define the meaning of the inequality \eqref{eq:intro-strict-cone} fails analytically. 
\begin{defn}
    Let $\cK_J^{an}(\alpha,\beta)$ consist of currents $T=S+\epsilon \chi$, for some closed positive current $S$ with analytic singularity,  K\"ahler form $\chi\in \beta$, and $\varepsilon>0$ such that 
\[[S]=\alpha-\epsilon\beta, \quad P_\chi(S)\leq c_{\alpha,\beta}.
    \]
  Here $S$ has analytic singularities means  locally \[
    S=\ddc u, \quad u=a\log\left(\sum_{\ell}|f_\ell|^2\right)+g,
 \qquad a>0,\quad f_\ell \text{ holomorphic, } g\in C^\infty.
    \]
\end{defn}

In the above definition of $\cK_J^{\rm an}(\alpha,\beta)$, we shift the class by $-\varepsilon\chi$. This is motivated by the mass concentration technique for $J$-equation developed by Chen \cite{C21}, which says that for any boundary classes $(\alpha,\beta)$, there is such a current $S\in \alpha-\varepsilon\beta$ with $P_\chi
(S)\leq \cab.$ Closed positive
$(1,1)$-currents can be approximated by currents with analytic singularity, with quantitative control of
their lower bounds.  However, the cone condition \eqref{eq:intro-strict-cone} is a
nonlinear inequality involving products of the regularized form, and
Demailly's regularization theorem does not by itself preserve it. That's the main difficulty for this problem. 

Set \[
\JEnK(\alpha,\beta)
 :=\bigcap_{T\in\mathcal K^{\rm an}_J(\alpha,\beta)}E_+(T).
\]
Here $E_+(T)=\{x\in X:\nu(T,x)>0\}$, where $\nu(T,x)$ is the Lelong number.  
It's not immediate from the definition that $\cK_{J}^{\rm an}(\alpha,\beta)$ is nonempty. If it is empty, we take $\JEnK=X$.
The main theorem of this paper is
\begin{thm}\label{thm:main-thm}Suppose $\dim_\C X=3$ and $(\alpha,\beta)$ is $J$-nef. Then:\begin{enumerate}
    \item $\JNull(\alpha,\beta)$ is a finite union
      $C_1\cup...\cup C_l\cup D_1\cup\cdots\cup D_m$ of irreducible curves and divisors;
\item there are $\varepsilon>0$ and
      $T=S+\varepsilon\chi\in\mathcal K^{\rm an}_J(\alpha,\beta)$
      such that
      \[
        E_+(T)=\JNull(\alpha,\beta).
      \]
\end{enumerate}
consequently
      \[\ \JEnK(\alpha,\beta)=\JNull(\alpha,\beta)\ .\]
\end{thm}
As a corollary, we prove the local smooth convergence of the $J$-flow under the $J$- nef assumption. The $J$-flow is defined as \begin{equation}
    \label{eq:J-flow}\dot\varphi=\cab-\tr_{\omega_\varphi}\chi, \qquad \varphi(0)=\varphi_0,
\end{equation}for some smooth $\varphi_0\in \Psh(X,\omega).$ Chen proved long-time existence of the flow \cite{C04}.  Weinkove
established convergence first on surfaces and then in higher dimension
under cohomological positivity assumptions
\cite{W04,W06}.  Song--Weinkove subsequently
identified the sharp smooth subsolution condition, proved convergence
to a critical metric under that condition, and analyzed the formation
of singularities when it fails \cite{SW08}.  The flow picture
was extended to inverse $\sigma_k$-equations by Fang--Lai--Ma and
Fang--Lai \cite{FLM11,FL12}; symmetry-reduced singular limits
were studied by Fang--Lai \cite{FL13}.  Collins--Sz\'ekelyhidi
proved that convergence depends only on the K\"ahler classes and gave a
numerical characterization on toric manifolds
\cite{CS17}. It is useful
to distinguish two different meanings of
``boundary'' that occur in the literature.  First, both background
classes remain K\"ahler, but the pair lies on the boundary of the
strict $J$-cone, so that the inequalities in
\eqref{eq:intro-strict-numerical} are nonnegative rather than positive.
On a surface this is equivalent to the class
$\rho=c\alpha-\beta$ being nef and big.  Fang--Lai--Song--Weinkove
treated the additional case in which $\rho$ has a smooth
semipositive representative: they obtained a uniform $C^0$ estimate
and smooth convergence of the $J$-flow, away from finitely many
negative curves, to a singular K\"ahler metric
\cite{FLSW14}.  Murakami later removed this semipositivity
assumption on compact K\"ahler surfaces.  For every 
$J$-nef pair he proved convergence in the sense of currents to the
unique current $T_\infty$ satisfying
\[
 cT_\infty-\chi\geq0,
 \qquad
 \left\langle(cT_\infty-\chi)^2\right\rangle=\chi^2
\]
in our notation \cite{M26a}.  Under additional geometric
semipositivity hypotheses, higher-dimensional boundary behavior was
studied by Sun \cite{S24}.  Very recently, Murakami proved a
boundary weak-solution and current-convergence theorem for generalized
Monge--Amp\`ere equations and their mixed Hessian flows; the
$J$-equation and $J$-flow are special cases
\cite{M26b}.  These are weak convergence theorems: by themselves
they do not identify the positive-Lelong locus of the limit or provide higher order partial regularity of the weak solution.

A second, logically different, use of the word ``degenerate'' allows
the fixed form appearing in the flow itself to be merely semipositive,
rather than K\"ahler.  Song--Weinkove constructed and proved convergence
of such a degenerate $J$-flow on minimal surfaces of general type
\cite{SW08}.  T\^o extended this to degenerate twisted
$J$-flows in arbitrary dimension, under a divisor-controlled
nondegeneracy assumption and a strict subsolution condition
\cite{T23}.  This degeneration of the background form should
not be confused with the $J$-nef boundary considered in
the present paper, where both $\alpha$ and $\beta$ remain K\"ahler.

\begin{cor}\label{cor:convergence-J-flow}Under the assumption of theorem \ref{thm:main-thm}, let $\varphi(t)$ solve the $J$-flow \eqref{eq:J-flow}. Set \[\overline{\varphi}(t):=\frac{\int_X\varphi(t)\chi^3}{\int_X\chi^3}.
\]Then there is a $\omega$-plurisubharmonic function $\varphi_\infty$, smooth on $X\setminus \JNull(\alpha,\beta)$, such that \[
\varphi(t)-\overline{\varphi}(t)\to \varphi_{\infty}\qquad \text{ in }C^\infty_{\rm loc}(X\setminus \JNull(\alpha,\beta)) \quad \text{ as }t\to \infty
.\]
Moreover, $\omega_{\varphi(t)}$ converges in the sense of currents to $\omega+\ddc\varphi_\infty.$
    
\end{cor}

We now give the rough ideas and motivations for the proof.  $J$-nefness is a boundary condition: all the required intersection numbers are nonnegative, but they vanish on the $J$-null locus.  We therefore perturb the cohomology class slightly into
the strict region, solve a smooth equation there, and then add positive
divisor currents to return to the original class.  These added currents
produce singularities exactly where the original numerical equalities
occur.  Schematically, the proof is
\[
 \begin{aligned}
 \text{boundary numerical data}
 &\longrightarrow \text{strict perturbed data}\\
 &\longrightarrow \text{smooth strict subsolution on a resolution}
 \longrightarrow \text{singular current}.
 \end{aligned}
\]
Null divisors provide perturbation directions directly on $X$.  A
null curve does not define a $(1,1)$-class, so we first resolve it and
use exceptional divisors above it.  The main difficulty is to find one
perturbation direction which makes all the vanishing inequalities
strict at the same time.

\smallskip
\noindent\emph{Step 1: perturb the null divisors.}
First suppose that there are no null curves.    Let $D_1,\ldots,D_m$ be the null divisors.  We look for an
effective $\mathbb Q$-divisor $ E=\sum_i a_iD_i, \text{ with }a_i>0$, 
and, set \[\alpha_s=\alpha-s\{E\}, \qquad r_s=\cab\alpha_s-\beta=:\rho-\cab s\{E\}.\]
For every null divisor $D_j$, a direct expansion gives
\begin{equation}\label{eq:intro-first-variation}
 J_{\cab}(D_j,\alpha_s,\beta)
 =-2s\,\rho ED_j+\cab s^2E^2D_j.
\end{equation}
Thus the equality on $D_j$ becomes strict for small $s>0$ if
$\rho ED_j<0$.  If $M_{ij}:=\rho D_iD_j,$
this means that we need a vector $a=(a_i)>0$ such that $Ma<0$.

\smallskip
\noindent\emph{Step 2: prove that the required direction exists.}
The key is the negative definiteness of $M$.  Its off-diagonal entries are nonnegative.  To prove first that $M\leq0$, for any effective integral divisor $F=\sum_ib_iD_i$, we construct a ray $\varphi_F(t):[0,\infty)\to \Psh(X,\omega)$ in the space of K\"ahler metrics in class $\alpha$ and relate the asymptotic slope of the $\cJ$-functional along this ray to the quadratic form $\sum_{i,j}M_{ij}b_ib_j$.
 Secondly, if $M$ had a kernel, the
same sign pattern would give an effective kernel divisor $E_0$.
We use the Hodge index theorem and Demailly-P\u{a}un theorem to show that $E_0$ is nef, contradicting
the key fact that $(\cab\alpha^2-2\alpha\beta)\cdot E_0=0$ and the strict positivity of $\cab\alpha^2-2\alpha\beta$ on every nonzero
modified-nef class.  Hence $M<0$.  Solving
$(-M)a=\mathbf1$, and then making a small rational approximation,
produces the required effective divisor $E$ in \emph{step 1}.

\smallskip
\noindent\emph{Step 3: solve the strict problem.}
For small $s>0$, both $\alpha_s$ and $r_s$ are K\"ahler.  We show that
\[
 (r_s^2-\beta^2)\cdot D=\cab J_{\cab}(D,\alpha_s,\beta)>0
\]
for every prime divisor $D$.  With $f_s={(r_s^3-3r_s\beta^2)}/{\beta^3}>0,$
Fang--Ma's theorem \cite{FM24} gives a K\"ahler form $R_s\in r_s$
such that
\[
 R_s^3=3R_s\wedge\chi^2+f_s\chi^3,
 \qquad R_s^2-\chi^2>0.
\]
It follows that
\[
 \omega_s=\frac{R_s+\chi}{c}\in\alpha_s,
 \qquad P_\chi(\omega_s)<c.
\]
After inserting the small shift $r_s-c\varepsilon\beta$, the same
argument gives
$\omega_{s,\varepsilon}\in
\alpha-s\{E\}-\varepsilon\beta$.  We now put back exactly the class
which was subtracted:
\[
 S=\omega_{s,\varepsilon}+s[E],\qquad
 T=S+\varepsilon\chi\in\mathcal K_J^{\rm an}(\alpha,\beta),
\]
where $[E]$ is the current of integration along $E$.  Adding this
positive current preserves the weak cone condition and creates the
desired analytic singularities along the null divisors.

\smallskip
\noindent\emph{Step 4: include the null curves.}
When null curves are present, $\rho$ is only nef and big, so the
preceding construction cannot be carried out directly on $X$.  We
resolve their union by $\mu:Y\to X$ and choose an exceptional divisor
$F$ such that $\mu^*\rho-s\{F\}$ is K\"ahler with an order-$s$
margin.  For every prime divisor $\Gamma\subset Y$, the projection
formula gives
\[
 \mu^*\rho\cdot\{F\}\cdot\{\Gamma\}=0.
\]
Thus the exceptional perturbation makes the class K\"ahler at first
order, but it does not yet cross the remaining divisor walls.  Those
walls are crossed at second order by adding
$-A_0s^2\mu^*\{E\}$, whose sign is controlled by $Ma<0$.  We can
then apply Fang--Ma's theorem  on $Y$.  The pushdown is smooth away from the
null curves; near them, local logarithmic models and Richberg
regularized maxima, in the spirit of Collins--Tosatti
\cite{CT15}, give neat analytic singularities while
preserving the weak cone inequality.  Finally we add the null-divisor
currents.  The resulting current lies in
$\mathcal K_J^{\rm an}(\alpha,\beta)$ and its positive-Lelong locus is
exactly the full $J$-null locus.  The reverse inclusion follows by
restriction and gluing near each null component.

The paper is organized according to these steps.  We first review the needed
notions, especially  Boucksom's divisorial Zariski decomposition.  We then prove one-side inclusion $\JNull\subset \JEnK$ in section 3. Next, we show the strict positivity of the evaluation of $\cab\alpha^2-2\alpha\beta$ on modified nef classes and deduce finiteness of irreducible $J$-null subvarieties in section 4.  In section 5 and 6, we prove the existence of divisorial direction $E$ where $J_{\cab}$ strictly increases. In section 7, we construct an element $T=S+\varepsilon\chi\in \cK_J^{\rm an}$ with $E_+(T)=\JNull(\alpha,\beta)$, and therefore close the proof of the main theorem. Finally, in section 8, we show the convergence of the $J$-flow away from $\JNull(\alpha,\beta).$

\section{Preliminaries}
We give some basic notations and facts that will appear in the proof. We fix a background K\"ahler metric $h$ on $X$.  \begin{itemize}
    \item Let $\eta\in H^{1,1}(X,\R)$. \begin{itemize}
        \item $\eta$ is a K\"ahler class if it contains a smooth K\"ahler form. K\"ahler class forms a cone $\cK_X$.
        \item $\eta$ is \emph{nef} if it lies in $\overline{\cK_X}$. Equivalently, for every $\epsilon>0$, $\eta$ contains a smooth closed form $\theta_\epsilon\geq -\epsilon h$
        \item $\eta$ is \emph{pseudo-effective} if it contains a closed positive
        current.  Such classes form the closed convex cone
        $\Psef(X)$.
        \item $\eta$ is \emph{big} if it contains a \emph{K\"ahler current},
        namely a closed current $T\in\eta$ such that
      $T\geq\delta h$ for some $\delta>0$.  The big cone is the
      interior of $\Psef(X)$.
      \item $\eta$ is \emph{modified K\"ahler} if it contains a K\"ahler current $T$
      such that the generic Lelong number $\nu(T,D)=0$ for every irreducible divisor $D$. Such class forms a cone $\MK(X)$.
\item $\eta$ is \emph{modified nef} if, for every $\epsilon>0$, it contains a
      closed current $T_\epsilon$ such that
      \[
        T_\epsilon\geq-\epsilon h,
        \qquad \nu(T_\epsilon,D)=0
        \quad\text{for every irreducible divisor }D.
      \]The cone of \emph{modified nef} class is denoted by $\MN(X)$.
    \end{itemize}
 \item Let $T$ be a closed $(1,1)$-current with $T\geq -ah$ for some constnat $a>0$. Demailly's regularization theorem \cite{D92} says that , for every $c>0$, there are currents $T_{c,k}\in \{T\}$, converges weakly to $T$, such that \begin{align*}
 T_{c,k}&\geq-\bigl(a+C c+\epsilon_k\bigr)h,
       \qquad \epsilon_k\downarrow0,\\
 \nu(T_{c,k},x)&=\max\{\nu(T,x)-c,0\},
\end{align*}
and $T_{c,k}$ is smooth outside the analytic set $
 E_c(T):=\{x:\nu(T,x)\geq c\}$. Here the constant $C$ depends  on the curvature of $(X,h)$.

There is also a direct one-line implication that will be used repeatedly.
If $\gamma\in\MN(X)$ and $\delta>0$, choose from the definition a current
\[
 T_{\delta/2}\in\gamma,\qquad
 T_{\delta/2}\geq-\frac{\delta}{2}h,\qquad
 \nu(T_{\delta/2},D)=0\quad\forall \text{ divisor }D.
\]
Then $
 T_{\delta/2}+\delta h\geq\frac{\delta}{2}h$
is a K\"ahler current in $\gamma+\delta h$, with the same generic
divisorial Lelong numbers.  Therefore, directly from the definitions,
\begin{equation}\label{eq:add-kahler}
 \gamma+\delta h\in\MK(X).
\end{equation}
Applying Demailly's regularization to $T_{\delta/2}+\delta h$ additionally
gives K\"ahler currents smooth away from subvarieties of codimension at least two.
\item A \emph{modification} is a proper surjective holomorphic map
\[
 \mu:\widehat X\longrightarrow X
\]
which is biholomorphic away from proper analytic subsets. 
Compositions of blow-ups are the basic examples(we will only use this in the paper). We say $E\subset \widehat X$ is an $\mu$-exceptional divisor if $\text{codim}_\C \mu(E)\geq 2.$
\item The following is a characterization of modified K\"ahler class by Boucksom \cite[proposition 2.3]{B04}:
A class $\eta$ is modified K\"ahler if and only if there exist a
modification $\mu:\widehat X\to X$ and a K\"ahler class
$\widehat\eta$ on $\widehat X$ such that
\[
 \eta=\mu_*\widehat\eta.
\]

Here $\mu_*\widehat\eta$ means: choose a K\"ahler form
$\widehat\theta\in\widehat\eta$, push it forward as a current, and take
the resulting cohomology class on $X$.  Although
$\mu_*\widehat\theta$ need not be smooth at the exceptional image, it is
a K\"ahler current and has zero generic Lelong number on every divisor,
because the singular-value locus of a modification has codimension at
least two.Combining with \eqref{eq:add-kahler}, we get exactly the form used later:
\begin{equation}\label{eq:boucksom-used}
 \gamma\in\MN(X),\ \delta>0
 \quad\Longrightarrow\quad
 \gamma+\delta h=\mu_{\delta *}\widehat h_\delta
\end{equation}
for a modification $\mu_\delta:Y_\delta\to X$ and a K\"ahler class
$\widehat h_\delta$ on $Y_\delta$.
 \item Boucksom's divisorial Zariski decomposition: Let $\xi\in\Psef(X)$.  Its generic minimal multiplicities define
\[
 N(\xi):=\sum_D\nu(\xi,D)\,D,\qquad
 Z(\xi):=\xi-\{N(\xi)\}.
\]
Boucksom calls $
 \xi=Z(\xi)+\{N(\xi)\}$, 
the \emph{divisorial Zariski decomposition}
\cite[Definition~3.7]{B04}. Here are some basic facts about the decomposition: \begin {enumerate}
    \item $Z(\xi)$ is \emph{modified nef} \cite[Proposition 3.8(i)]{B04}.
    \item $N(\xi)$ is an effective $\R$-divisor with only finitely many
      components
      \cite[Theorem~3.12(i)]{B04};
      \item those components form an \emph{exceptional family}: their positive
      span meets $\MN(X)$ only at $0$
      \cite[Definition~3.10 and Proposition~3.11]{B04};
\item the classes of the prime divisors in an exceptional family are
      linearly independent in $\NS(X)_\R$
      \cite[Proposition~3.11(iii)]{B04}.
\end{enumerate}
Here $\NS(X)_\R$ is the real span in $H^{1,1}_{BC}(X,\R)$ of first
Chern classes of holomorphic line bundles.  Every divisor class belongs
to this space.  Its dimension
$\rho(X)$ is the \emph{Picard number}; it is finite.  Therefore an exceptional family has at most
$\rho(X)$ members.  This simple finite-dimensional fact is what makes
the $J$-null locus finite in lemma \ref{lem:finite-null-divisor}.


\end{itemize}

\section{Null locus forces singularities}

In this section, we prove an easy direction of the main theorem: $\JNull(\alpha,\beta)\subset \JEnK(\alpha,
\beta).$ We start with the the strict cone metric in a perturbed class. 
\begin{lemma}
\label{lem:perturbed-subsolutions}
For every $t>0$, there is a smooth K\"ahler form
\[
 \Omega_t\in(1+t)\alpha
\]
such that
\[
 c_{\alpha,\beta}\Omega_t^2-2\chi\wedge\Omega_t>0,
 \qquad\text{equivalently}\qquad P_\chi(\Omega_t)<c_{\alpha,\beta}.
\]
\end{lemma}
\begin{proof}
    For $t>0$, set $A_t=(1+t)\alpha$ and $
 \delta_t:=\frac{\cab t}{2(1+t)}.$ 
For an irreducible $p$-dimension subvariety $V$, $1\leq p\leq3$, one has
\begin{align*}
 &\bigl(\cab-(3-p)\delta_t\bigr)A_t^p
       -p\beta A_t^{p-1}\\
 &\quad=(1+t)^{p-1}
 \left[
   \cab\alpha^p-p\beta\alpha^{p-1}
   +\bigl(\cab t-(3-p)\delta_t(1+t)\bigr)\alpha^p
 \right].
\end{align*}
The last coefficient is $0$, $\cab t/2$, and $\cab t$ for
$p=1,2,3$, respectively.  Thus all of Chen's uniform numerical
inequalities hold on every irreducible  $p$-dimension subvariety.  

Choose the positive constant $f_t$ by $
 \frac{f_t}{6}\beta^3
 =\frac{ c_{\alpha,\beta}}{6} A_t^3-\frac12\beta A_t^2.$ 
The normalization $c_{\alpha,\beta}\alpha^3=3\beta\alpha^2$ gives
\[
 \frac{c_{\alpha,\beta}}{6} A_t^3-\frac12\beta A_t^2
 =\frac{c\alpha^3}{6}\,t(1+t)^2>0.
\]
Thus $f_t$ is smooth and satisfies both remaining hypotheses of
Chen's Theorem~1.11 \cite{C21}:
\[
 f_t>0>-\frac16c^{-2},
 \qquad
 \int_X f_t\frac{\chi^3}{3!}
 =c\int_X\frac{A_t^3}{3!}
   -\int_X\chi\wedge\frac{A_t^2}{2!}>0.
\]
That theorem therefore produces a K\"ahler form
$\Omega_t\in A_t$ satisfying
\[
 c_{\alpha,\beta}\Omega_t^2-2\chi\wedge\Omega_t>0,\qquad
 \operatorname{tr}_{\Omega_t}\chi
      +f_t\frac{\chi^3}{\Omega_t^3}=c_{\alpha,\beta},
\]
which proves the lemma.

\end{proof}

\begin{prop}\label{prop:easy-direction}
Every
$T\in\mathcal K^{\rm an}_J(\alpha,\beta)$ satisfies
\[
 \JNull(\alpha,\beta)\subset E_+(T).
\]
\end{prop}
\begin{proof}
Write $T=S+\varepsilon\chi$ and $Z=E_+(T)$.  On
$X\setminus Z$, the currents $S,T$ are smooth.  If
$\mu_1,\mu_2$ are the eigenvalues of $\chi$ relative to $S$ on a
complex two-plane, then $\mu_1+\mu_2\leq c_{\alpha,\beta}.$ 
The corresponding eigenvalues relative to
$T=S+\varepsilon\chi$ are
$\mu_i/(1+\varepsilon\mu_i)$.  Hence
\begin{equation}\label{eq:margin}
 P_\chi(T)\leq \frac{c_{\alpha,\beta}}{1+\epsilon c_{\alpha,\beta}/2}=:c_\epsilon<c,\qquad T\geq\varepsilon\chi.
\end{equation}

Let $V$ be a $J$-null irreducible subvariety of dimension
$p=1$ or $2$, and suppose that $V\not\subset Z$.  Choose
coordinate balls
\[
 U'\Subset U\Subset X\setminus Z
\]
such that $U'\cap V_{\rm reg}\neq\varnothing$.
For $t>0$, let $\Omega_t\in(1+t)\alpha$ be the smooth strict
subsolution from Lemma~\ref{lem:perturbed-subsolutions}.  Since
$(1+t)T$ and $\Omega_t$ are cohomologous, on $X$ we may write
\[
 (1+t)T=\Omega_t+\ddc\varphi_t,
\]
where $\varphi_t$ is smooth on $X\setminus Z$ and tends to
$-\infty$ along $Z$ when $Z\neq\varnothing$.  If
$Z=\varnothing$, take $\widetilde\Omega_t=(1+t)T$.  If
$Z\neq\varnothing$, proceed as follows.
Fix a small regularization width $\tau>0$.  After adding a constant
to $\varphi_t$, we can
arrange
\[
 \varphi_t>2\tau\quad\hbox{on }U',
 \qquad
 \varphi_t<-2\tau\quad\hbox{on a neighborhood of }Z.
\]
Indeed, $\varphi_t$ is bounded on $\overline{U'}$ and tends
uniformly to $-\infty$ near the compact analytic set $Z$.  Moreover,
\[
 P_\chi((1+t)T)=\frac1{1+t}P_\chi(T)
 \leq\frac{c_\varepsilon}{1+t}<c, \text{ on } X\setminus Z.
\]
    Define $\tilde{\varphi}_t:=\max_{\tau}\{0,\varphi_t\}$, where $\max_\tau$ is the regularized maximum function. We have $\tilde{\varphi}_t=\varphi_t$ on $U'$, and $\tilde{\varphi}_t=0$ in a neighborhood of $Z$. By the convexity of the $P_\chi(\cdot)$, we get $\tilde{\Omega}_t=\Omega_t+\ddc\tilde{\varphi}_t\in (1+t)\alpha$ with $P_\chi(\tilde{\Omega}_t)<c_{\alpha,\beta}$, and on $U'$, $\tilde{\Omega}_t=(1+t)T$.
    The form $
 c_{\alpha,\beta}\widetilde\Omega_t^p
   -p\chi\wedge\widetilde\Omega_t^{p-1}$ 
is positive on $V_{\rm reg}$.  On the fixed open set
$U'\cap V_{\rm reg}$, the equality
$\widetilde\Omega_t=(1+t)T$, the bound
$P_\chi(T)\leq c_\epsilon$, and
$T\geq\varepsilon\chi$ give the uniform lower estimates
\[
\begin{aligned}
 c_{\alpha,\beta}\widetilde\Omega_t-\chi
 &\geq(c_{\alpha,\beta}-c_\varepsilon)\varepsilon\chi
 &&(p=1),\\
 c_{\alpha,\beta}\widetilde\Omega_t^2
     -2\chi\wedge\widetilde\Omega_t
 &\geq(c_{\alpha,\beta}-c_\varepsilon)\varepsilon^2\chi^2
 &&(p=2).
\end{aligned}
\]
Indeed, $P_\chi(T)\leq c_\varepsilon$ implies both
\[
 \chi\leq c_\varepsilon T,\qquad
 c_\varepsilon T^2-2\chi\wedge T\geq0.
\]
Writing $\lambda=1+t$, the first lower bound follows from
\[
 c\lambda T-\chi
 =(c-c_\varepsilon+ct)T+(c_\varepsilon T-\chi)
 \geq(c-c_\varepsilon)T,
\]
and the second from
\[
\begin{aligned}
 c\lambda^2T^2-2\lambda\chi\wedge T
 &=\lambda\bigl[
   (c_{\alpha,\beta}-c_\varepsilon+c_{\alpha,\beta}t)T^2
   +(c_\varepsilon T^2-2\chi\wedge T)\bigr]\\
 &\geq(c_{\alpha,\beta}-c_\varepsilon)T^2.
\end{aligned}
\]
Finally,
\[
 T^2-\varepsilon^2\chi^2
 =(T-\varepsilon\chi)\wedge(T+\varepsilon\chi)\geq0.
\]

Set
\[
 \delta_0:=
 \begin{cases}
 (c_{\alpha,\beta}-c_\varepsilon)\varepsilon\displaystyle\int_{U'}\chi,
       &p=1,\\[3mm]
 (c_{\alpha,\beta}-c_\varepsilon)\varepsilon^2\displaystyle\int_{U'}\chi^2,
       &p=2.
 \end{cases}
\]
Then $\delta_0>0$ and is independent of $t$. 
 Integration on $V$ therefore gives 
\begin{equation}\label{eq:positive-lower-bound}
 J_{c_{\alpha,\beta}}(V,(1+t)\alpha,\beta)
 =\int_V\bigl(c_{\alpha,\beta}\widetilde\Omega_t^p
     -p\chi\wedge\widetilde\Omega_t^{p-1}\bigr)\geq \delta_0. 
\end{equation}
On the other hand, $J_{c_{\alpha,\beta}}(V,\alpha,\beta)=0$, so a direct
cohomological calculation gives
\[
\begin{aligned}
 J_{c_{\alpha,\beta}}(V,(1+t)\alpha,\beta)
 &=(1+t)^{p-1}
   \bigl(c_{\alpha,\beta}(1+t)\alpha^p-p\beta\alpha^{p-1}\bigr)\cdot V\\
 &=c_{\alpha,\beta}t(1+t)^{p-1}\alpha^p\cdot V
 \longrightarrow0.
\end{aligned}
\]
This contradicts \eqref{eq:positive-lower-bound}.  Hence $V\subset Z$.
\end{proof}
\section{Finiteness of the null locus}
Set \[
\rho=c_{\alpha,\beta}\alpha-\beta, \qquad z:=\rho^2-\beta^2=\cab(\cab\alpha^2-2\alpha\beta).\] Then direct computation using $\cab\alpha^3=3\alpha\beta^2$ gives  
 \[\rho^3-3\rho\beta^2=2\beta^3>0.\]
We first show that the class $\rho$ is nef and big, then the finiteness of the null curves follows from Collins-Tosatti's theorem.
\begin{lemma}\label{lem:rho-nef}
Assume $(\alpha,\beta)$ is $J$-nef.  Then $\rho=\cab\alpha-\beta$ is nef and
big.  Moreover,
\[
 \operatorname{Null}(\rho)
 =C_1\cup\cdots\cup C_\ell
\]
is a finite union of irreducible curves, and these are exactly the
$J$-null curves.
    
\end{lemma}
\begin{proof}
For $s>0$, set
\[
 \rho_s:=\rho+s\alpha.
\]
If $C$ is an irreducible curve,  $J$-nefness gives
\[
 \rho_s\cdot C=\rho\cdot C+s\alpha\cdot C>0.
\]
If $D$ is a prime divisor, then
\[
 \rho^2\cdot D\geq\beta^2\cdot D>0
\]
because $z\cdot D\geq0$.  The divisor inequality also gives
\[
 \rho\alpha\cdot D
 =c_{\alpha,\beta}\alpha^2\cdot D-\beta\alpha\cdot D
 \geq\frac{c_{\alpha,\beta}}{2}\alpha^2\cdot D>0.
\]
Consequently,
\[
 \rho_s^2\cdot D
 =\rho^2\cdot D+2s\rho\alpha\cdot D+s^2\alpha^2\cdot D>0.
\]
It remains to check the threefold $X$. By Khovanskii--Teissier inequality \cite[proposition 2.1]{DP03}, one has \[
(\beta\alpha^2)^2\geq(\alpha^3)(\beta^2\alpha), \quad (\beta^2\alpha)^2\geq (\beta\alpha^2)(\beta^3).
\]
It follows that $(\beta\alpha^2)(\beta^2\alpha)\geq (\alpha^3)(\beta^3).$ Hence \[
\rho\beta^2=\frac{3(\beta\alpha^2)(\beta^2\alpha)}{\alpha^3}-\beta^3\geq 2\beta^3>0.
\]
One has $\rho^3=3\rho\beta^2+2\beta^3>0$.\
Hence \[
\rho_s^3=\rho^3+3s\rho^2\alpha+3s^2\rho\alpha^2+s^3\alpha^3=\rho^3+s(\cab^2\alpha^3+3\beta^2\alpha)+2s^2\cab\alpha^3+s^3\alpha^3>0.
\]
Thus $\rho_s$ has positive top self-intersection on every positive-dimensional irreducible analytic subvariety of $X$. For
$s\gg1$, the class $\rho_s$ is K\"ahler.  The connected ray
$\{\rho_s:s>0\}$ lies in the numerical positive cone of
Demailly--P\u{a}un, and it meets the K\"ahler cone.  Their
connected-component theorem therefore shows that $\rho_s$ is
K\"ahler for every $s>0$\cite{DP04}.  Letting $s\downarrow0$ proves that
$\rho$ is nef.  Since $\rho^3>0$, the nef-and-big criterion of
Demailly--P\u{a}un shows that $\rho$ is big. 

Collins--Tosatti's null-locus theorem now gives
\[
 E_{\rm nK}(\rho)
 =\operatorname{Null}(\rho)
 =\bigcup_{\rho^{\dim V}\cdot V=0}V.
\]
There is no divisorial component, because
\[
 \rho^2\cdot D\geq\beta^2\cdot D>0
\]
for every prime divisor $D$.  Hence the positive-dimensional
irreducible components are finitely many curves $C_j$.  For a curve,
$\rho\cdot C=J_{\cab}(C,\alpha,\beta)$, so they are exactly the
$J$-null curves.

\end{proof}

\begin{lemma}[Strictly positivity on modified-nef classes]\label{lem:MN}
For every nonzero $\gamma\in\MN(X)$,
\[
 z\cdot \gamma=\cab(c_{\alpha,\beta}\alpha^2-2\beta\alpha)\cdot\gamma>0.
\]
\end{lemma}

\begin{proof}
Choose $s>1$ sufficiently close to $1$ that $
 \rho^3-3s^2\rho\beta^2>0$; 
this is possible since $\rho^3-3\rho\beta^2=2\beta^2>0$ and the left-hand side is continuous in $s$. Fix a K\"ahler class $h$ and $\delta>0$.  By
\eqref{eq:add-kahler}, $\gamma+\delta h$ is modified K\"ahler.
Boucksom's characterization, therefore, gives a modification
$\mu_\delta:Y_\delta\to X$ and a K\"ahler class
$\widehat h_\delta$ on $Y_\delta$ such that
\[
 (\mu_\delta)_*\widehat h_\delta=\gamma+\delta h.
\]

Set
\[
 A_\delta:=\mu_\delta^*\rho,\qquad
 B_\delta:=s\mu_\delta^*\beta.
\]
Pullback preserves nefness, so $A_\delta$ is nef because $\rho$ is
nef by lemma~\ref{lem:rho-nef}.  Moreover $B_\delta$ is nef:
it is represented by the pullback of the K\"ahler form $s\beta$,
which is smooth and semipositive.  Pullback preserves top intersection numbers,
and hence
\[
 A_\delta^3-3A_\delta\cdot B_\delta^2
 =\rho^3-3s^2\rho\beta^2>0.
\]
By \cite[lemma 2.4]{X16}, there is a closed strictly positive
current $\Theta_\delta$ whose class is
\[
 \{\Theta_\delta\}=A_\delta^2-B_\delta^2
 =\mu_\delta^*(\rho^2-s^2\beta^2).
\]
Because $\widehat h_\delta$ is K\"ahler and
$\Theta_\delta$ is strictly positive,
\[
 \int_{Y_\delta}\Theta_\delta\wedge\widehat h_\delta>0.
\]
The projection formula converts this analytic positivity upstairs into
the required intersection inequality downstairs:
\[
\begin{aligned}
 (\rho^2-s^2\beta^2)\cdot(\gamma+\delta h)>0.
\end{aligned}
\]
Indeed, the left-hand side equals
\[
 \mu_\delta^*(\rho^2-s^2\beta^2)\cdot\widehat h_\delta
 =(A_\delta^2-B_\delta^2)\cdot\widehat h_\delta
 =\int_{Y_\delta}\Theta_\delta\wedge\widehat h_\delta.
\]
Letting $\delta\downarrow0$ gives
\[
 (\rho^2-s^2\beta^2)\cdot\gamma\geq0.
\]
A modified-nef class is pseudo-effective: take
$\epsilon_j\downarrow0$, choose
$T_j\in\gamma$ with $T_j\geq-\epsilon_jh$, normalize the potentials,
and take a weakly convergent subsequence.  Its limit is a positive current
in $\gamma$.  Since $\gamma\neq0$ and $\beta$ is K\"ahler,
 $\beta^2\cdot\gamma>0$.  Finally,
\[
 z\cdot\gamma
 =(\rho^2-s^2\beta^2)\cdot\gamma
   +(s^2-1)\beta^2\cdot\gamma>0.
\]
\end{proof}

Now we can prove the finiteness of null divisors. 
\begin{lemma}\label{lem:finite-null-divisor}
There are only finitely many $J$-null prime divisors. Denote them by $D_1,...,D_m$, and
\begin{equation}\label{eq:zeroface}
 \{\xi\in\Psef(X):(c_{\alpha,\beta}\alpha^2-2\beta\alpha)\cdot\xi=0\}
 =\sum_{i=1}^m\R_{\geq0}\{D_i\}.
\end{equation}
\end{lemma}
\begin{proof}
    For $\xi\in \Psef(X)$, take its divisorial Zariski decomposition $\xi=Z(\xi)+N(\xi).$ Then $Z(\xi)$ is \emph{modified nef}, and $N(\xi)=\sum_{\ell=1}^q a_\ell E_\ell$, with $a_\ell>0$. Lemma~\ref{lem:MN} says that $z\cdot Z(\xi)>0$ unless $Z(\xi)=0$.  J-nefness says that $z\cdot D\geq 0$ for every irreducible divisor $D$. In particular \[
    z\cdot \xi=z\cdot Z(\xi)+\sum_{\ell=1}^qb_\ell z\cdot E_\ell\geq 0.
    \]
    It follows that
$z\cdot\xi=0$ exactly when
\begin{equation}\label{eq:zero-zariski-parts}
 \begin{aligned}
  &Z(\xi)=0,\quad\text{and}\\
  &\operatorname{Supp}N(\xi)
    \text{ consists of J-null prime divisors}.
 \end{aligned}
\end{equation}
Any finite collection of J-null irreducible divisors is exceptional in
Boucksom's sense: its positive span meets $\MN(X)$ only at $0$, by
Lemma~\ref{lem:MN}.  Boucksom's Proposition~3.13 says that the classes in
an exceptional family are linearly independent.  Thus the cardinality of
every finite collection of null irreducible divisors is bounded by the Picard
number $\rho(X)=\dim\NS(X)_\R$.  If there were infinitely many null
prime divisors, choosing $\rho(X)+1$ distinct ones would contradict this
bound.  There are therefore only finitely many.  Calling them
$D_1,\ldots,D_m$, \eqref{eq:zero-zariski-parts} proves
\eqref{eq:zeroface}.
\end{proof}

\section{the first variation of $J$-slope in null-divisor direction}
For the J-Null irreducible divisor $D_1,...,D_m$ in the previous lemma, set \[M_{ij}:=\rho\cdot D_i\cdot D_j=\int_X(c\alpha-\beta)\wedge \theta_{D_i}\wedge\theta_{D_j},\]
where $\theta_D$ is the Chern curvature of the holomorphic line bundle $\cO_X(D)$ equipped with a hermitian metric. 
We first show that $M$ is negative definite on $\Z_{\geq 0}^m$.
\begin{lemma}\label{lem:seminegative-M}
    If $F=\sum_{i}m_iD_i$, with $0\leq m_i\in \Z$. Then $(c_{\alpha,\beta}\alpha-\beta)\cdot F^2\leq 0$
\end{lemma}
\begin{proof}
    The proof is divided into 4 steps. 
    
\noindent\emph{1. Approach the boundary by coercive $J$-functionals.} This time, we perturb the class $\beta$ and keep 
    For $t>0$, put \[
    \chi_t=\chi+t\omega, \qquad \beta_t:=\beta+t\alpha,\qquad c_t:=\cab+3t.
    \]The number $c_t$ is the normalized $J$-slope: $c_t=\frac{3\beta_t\alpha^2}{\alpha^3}$.
For a $p$-dimensional irreducoble subvariety, Chen's uniform numerical condition \cite{C21} is \[
c_t\alpha^p-p\beta_t\alpha^{p-1}=(3-p)t\alpha^p+c\alpha^p-p\beta\alpha^{p-1}\geq (3-p)t\alpha^p.
\]
Recall the definition of $\cJ$-functional: for a K\"ahler form $\gamma$, and $\varphi\in \Psh(X,\omega)$, \begin{equation}\label{eq:J-functional}
    \cJ_\gamma(\varphi):=\frac{1}{6}\sum_{j=0}^2\int_X\varphi\gamma\wedge\omega^{j}\wedge_\varphi^{2-j}-\frac{c_\gamma}{24}\sum_{j=0}^3\int_X\varphi\omega^j\wedge\omega_\varphi^{3-j}, \quad \text{ where }c_\gamma=\frac{3\gamma\cdot \alpha^2}{\alpha^3}.
\end{equation}
By Chen's theorem, $\cJ_{\chi_t}$ is proper, i.e.  \begin{equation}\label{eq:Chen-coercivity}
    \cJ_{\chi_t}(\varphi)\geq \delta_t\cJ_\omega(\varphi)-C_t, \qquad \text{ for some constant }\delta>0,C_t>0, \text{ and for all }\varphi\in \Psh(X,\omega). 
\end{equation}
\smallskip
\noindent\emph{2. Smooth approximation for the logarithmic singularity.}
 Fix $\lambda>0$ sufficiently small, such that $\alpha-\lambda[F]$ is a K\"ahler class. We choose a K\"ahler form $\omega_\lambda$ in that class. Then the form $\theta_F:=\lambda^{-1}(\omega-\omega_\lambda)$ represents $c_1(\cO_X(F))$. Let $s_F$ be a defining section of $\cO_X(F)$, and choose the hermitian metric $h$ on it so that $|s_F|_h^2\leq 1$ and the Chern curvature of $h$ is exactly $\theta_F$, i.e. on $X\setminus F$, \[
-\ddc\log(|s_F|_h^2)=\theta_F=\lambda^{-1}(\omega-\omega_\lambda),
\]  and globally as distribution $\ddc\log|s_F|_h^2=[F]-\theta_F$. For $R>0$, set \[
u_R:=\log(e^{-R}+|s_F|_h^2), \qquad \varphi_R:=\lambda u_R.
\]
Then direct computation shows that\begin{equation}
    \ddc u_R=\frac{|s_F|_h^2}{e^{-R}+|s_F|_h^2}\ddc\log(|s_F|_h^2)-d\left(\frac{e^{-R}}{e^{-R}+|s_F|_h^2}\right)\wedge d^c\log(|s_F|_h^2).
\end{equation} 
We set $a_R:=\frac{e^{-R}}{e^{-R}+|s_F|_h^2}$, and denote the second term on the right-hand side of the equation for $\ddc u$ above by $S_R$ Then $S_R=-da_R\wedge d^c\log|s_F|_h^2$, which extends smoothly over $F$, is semipositive on $X$, and has rank at most $1$. Hence $S_R^2=0$. It follows that \begin{equation}
    \omega_{\varphi_R}=\omega+\ddc\varphi_R=a_R\omega+(1-a_R)\omega_\lambda+\lambda S_R>0.
\end{equation}
Thus $\varphi_R\in \Psh(X,\omega)$.

\smallskip
\noindent\emph{3. The local logarithmic-slope calculation.} We set \[
T_R:=\theta_F+\ddc u_R=a_R\theta_F+
S_R.
\]
We claim that, for every smooth closed form $\Gamma$ of complementary
bidegree and $k=1,2,3$,
\begin{equation}\label{eq:basic-log-slope}
 \lim_{R\to\infty}\frac1R
  \int_Xu_R\,T_R^k\wedge\Gamma
 =-\int_X\theta_F^k\wedge\Gamma=:-F^k\cdot\{\Gamma\}.
\end{equation}
For $k=0$, the corresponding limit is zero because logarithmic
singularities are locally integrable and $u_R\to\log |s_F|_h^2$ in
$L^1_{\rm loc}$.

Write
\[
 I_k(R):=\int_Xu_R\,T_R^k\wedge\Gamma.
\]
Since $\dot u_R=-a_R$, differentiation followed by integration by
parts gives
\begin{equation}\label{eq:Ik-derivative}
\begin{aligned}
 \frac{dI_k}{dR}
 =&-\int_Xa_RT_R^k\wedge\Gamma+k\int_Xu_R(-\ddc a_R)T_R^{k-1}\wedge\Gamma\\
 =&-(k+1)\int_Xa_RT_R^k\wedge\Gamma+k\int_Xa_R\theta_F\wedge T_R^{k-1}\wedge\Gamma=I+II.
\end{aligned}
\end{equation}
It remains to take the limit in $I,II$. First, we prove that for any integer $j\geq 0$, \begin{equation}\label{eq:aP-current-limit}
    a_R^j S_R=-a_R^jda_R\wedge d^c\log|s_F|_h^2\wto \frac{1}{j+1}[F], \qquad \text{as }R\to \infty.  
\end{equation}
The intuition of the convergence is:  at a smooth point of $F$, move in the normal direction from the divisor to its completement. Then $a_R$ decreases from $1$ to $0$. So the radial part of the mass is simply $-\int_1^0a^jda=\frac{1}{j+1}$. The angular part $d^c\log|s_F|_h^2$ is the Poincare-Lelong mass of $F$. This gives the right side of \eqref{eq:aP-current-limit}. Now we prove it rigorously. Define \[
\Psi_{j,R}(s)
 :=\frac1{j+1}\int_1^s
   \frac{1-\bigl(e^{-R}/(e^{-R}
+t)\bigr)^{j+1}}{t}\,dt.
\]
The apparent singularity of the integrand at $s=0$ is removable, so
$\Psi_{j,R}(|s_F|_h^2)$ is smooth on $X$.  Direct differentiation gives\[\begin{aligned}
 \ddc\Psi_{j,R}(|s_F|_h^2)&=\frac{1}{j+1}d\left((1-a_R^{j+1})d^c\log|s_F|_h^2\right) =a_R^jS_R  +\frac{1-a_R^{j+1}}{j+1}\ddc\log|s_F|_h^2.
\end{aligned}
\]
It follows that \begin{equation}\label{eq:aRS_R-limit}
   a_R^jS_R=\ddc\Psi_{j,R}(|s_F|_h^2)+\frac{1-a_R^{j+1}}{j+1}\theta_F. 
\end{equation}

Note that $|s_F|_h^2\leq 1$, one has \[
|\Psi_{j,R}(|s_F|_h^2)|\leq \frac{1}{j+1}|\log|s_F|_h^2|.
\] Moreover, $\Psi_{j,R}(|s_F|_h^2)\to \frac{1}{j+1}\log |s_F|_h^2$ pointwise on $X\setminus F$.  Dominated convergence and Poincaré-Lelong formula imply \[
 dd^c \Psi_{j,R}(|s_F|_h^2)\wto \frac{1}{j+1}([F]-\theta_F).
 \]
 The second term in \eqref{eq:aRS_R-limit} converges in $L^1$ to $\frac{1}{j+1}\theta_F$, which proves \eqref{eq:aP-current-limit}.

Since $S_R^2=0$, we now have \[
\begin{aligned}
 a_RT_R^k
 &=a_R^{k+1}\theta_F^k
   +k\,a_R^kS_R\theta_F^{k-1},\\
 a_R\theta_FT_R^{k-1}
 &=a_R^k\theta_F^k
   +(k-1)a_R^{k-1}S_R\theta_F^{k-1}.
\end{aligned}
\]
The terms containing only powers of $a_R$ tend to zero in $L^1$,
whereas \eqref{eq:aP-current-limit} evaluates the other terms.
Consequently, as currents,
\begin{equation}\label{eq:weighted-current-limits}
\begin{aligned}
 a_RT_R^k
 &\longrightarrow\frac{k}{k+1}
       [F]\wedge\theta_F^{k-1},\\
 a_R\theta_FT_R^{k-1}
 &\longrightarrow\frac{k-1}{k}
       [F]\wedge\theta_F^{k-1},
\end{aligned}
\end{equation}
where the second limit is zero when $k=1$.
Substitution of \eqref{eq:weighted-current-limits} into
\eqref{eq:Ik-derivative} gives
\[
 \lim_{R\to\infty}I_k'(R)
 =-F^k\cdot\{\Gamma\}.
\]
This proves \eqref{eq:basic-log-slope}.  

Now we consider $\omega_{\varphi_R}=\omega_\lambda+\lambda T_R$. By \eqref{eq:basic-log-slope}, we have \begin{equation}\label{eq:log-binomial-slope}
\begin{aligned}
 &\lim_{R\to\infty}\frac1R
   \int_X\lambda u_R\,\Gamma
       \wedge\omega_{\varphi_R}^{p}\\
       =&-\lambda\,\{\Gamma\}\cdot
   \sum_{k=1}^{p}\binom pk
       \lambda^kF^k(\alpha-\lambda F)^{p-k}.\\
       =&-\lambda\,\{\Gamma\}\cdot
   \bigl(\alpha^p-(\alpha-\lambda F)^p\bigr).
\end{aligned}
\end{equation}
For $p=0$, both sides are zero.

\smallskip
\noindent\emph{4. Compute the slope of the $\cJ$-functional.} Apply \eqref{eq:log-binomial-slope} term by term in \eqref{eq:J-functional}, we get \begin{equation}\label{eq:asymptotic-gamma-slope-polynomial}\begin{aligned}
    \lim_{R\to \infty}\frac{1}{R}\cJ_\gamma(\varphi_R)=&\frac{1}{6}\sum_{j=0}^1-\lambda\{\gamma\}\cdot\alpha^j\cdot(\alpha^{2-j}-(\alpha-\lambda F)^{2-j})\\
    &-\frac{c_{\gamma}}{24}\sum_{j=0}^2-\lambda\alpha^j\cdot(\alpha^{3-j}-(\alpha-\lambda F)^{3-j})\\
    =&\frac{1}{6}\lambda\{\gamma\}\cdot(-3\lambda\alpha\cdot F+\lambda^2F^2)+\frac{c_\gamma}{24}\lambda(6\lambda\alpha^2F-4\lambda^2\alpha F^2+\lambda^3F^3)\\
    =&(\frac{1}{4}c_\gamma\alpha^2-\frac{1}{2}\{\gamma\}\alpha)\lambda^2F-\frac{1}{6}(c_\gamma\alpha-\{\gamma\})\lambda^3F^2+\frac{c_\gamma}{24}\lambda^4F^3
\end{aligned}
\end{equation}
Similarly, one has\begin{equation}\label{eq:asymptotic-omega-slope-polynomial}
\lim_{R\rightarrow \infty}\frac{1}{R}\cJ_\omega(\varphi_R)=\frac{1}{4}\alpha^2\lambda^2F-\frac{1}{3}\alpha\lambda^3F^2+\frac{1}{8}\lambda^4F^3. 
\end{equation}
Note that $\alpha^2F=\int_{F}\alpha^2>0$. Hence for sufficient small $\lambda>0$, one has \[
\lim_{R\rightarrow \infty}\frac{1}{R}\cJ_\omega(\varphi_R)>0.
\]
 The coercivity of $\cJ_{\chi_t}$ \eqref{eq:Chen-coercivity}  gives 
 \[
 \lim_{R\to \infty}\frac{1}{R}\cJ_{\chi_t}(\varphi_R)\geq \delta_t\lim_{R\to \infty}\frac{1}{R}\cJ_\omega(\varphi_R)>0.
 \]
 Let $t\rightarrow 0$, we get \begin{equation}
     \label{eq:seminegative}
(\frac14c_{\alpha,\beta}\alpha^2-\frac12\alpha\beta)\lambda^2F
   -\frac16(c_{\alpha,\beta}\alpha-\beta) \lambda^3 F^2
   +\frac{c_{\alpha,\beta}}{24}\lambda^4F^3\ge0 \quad \text{ for all sufficient small }\lambda>0.
  \end{equation}
Since every component of $F$ is $J$-null, \[(
c_{\alpha,\beta}\alpha^2-2\alpha\beta)F=0.
\]
Thus, after division by $\lambda^3$, \eqref{eq:seminegative} becomes \[
-\frac{1}{6}\rho F^2+\frac{\cab}{24}\lambda F^3\geq 0.
\]
Letting $\lambda\to 0$ gives $\rho F^2\leq 0$. 
\end{proof}

Our goal is show that the matrix $M$ is  negative definite rather than only semi-negative definite. 
We will need two elementary lemmas about symmetric matrix. \begin{lemma}\label{lem:sign-matrix}
Let $A=(A_{ij})$ be a real symmetric matrix with
$A_{ij}\geq0$ whenever $i\neq j$.
\begin{enumerate}[label=\textup{(\roman*)}]
\item An eigenvector for the largest eigenvalue of $A$ can be chosen
      with all coordinates nonnegative.
\item If $M$ is irreducible.  Then an eigenvector
      for the largest eigenvalue can be chosen with all coordinates
      strictly positive.
\end{enumerate}
\end{lemma}
Here a matrix $A$ is called irreducible if no simultaneous permutation of the rows and columns
      writes $A$ in the form
      \[
        \begin{pmatrix}A^{(1)}&0\\0&A^{(2)}\end{pmatrix},
      \]
      where both square blocks have positive size
\begin{proof}
Let $\lambda_{\max}$ be the largest eigenvalue and choose a real
eigenvector $v\neq0$ for it.  Put $w_i=|v_i|$.  The diagonal terms
of $v^TAv$ and $w^TAw$ are identical, while for $i<j$
\[
 A_{ij}|v_i||v_j|\geq A_{ij}v_iv_j.
\]
Consequently
\[
 w^TAw\geq v^TAv
          =\lambda_{\max}\lVert v\rVert^2
          =\lambda_{\max}\lVert w\rVert^2.
\]
The reverse inequality follows from the maximal characterization of
$\lambda_{\max}$.  Equality therefore holds, and $w$ is also a
largest-eigenvalue eigenvector.  This proves (i).

For (ii), let
\[
 I_0:=\{i:w_i=0\},
 \qquad
 I_+:=\{i:w_i>0\}.
\]
The set $I_+$ is nonempty because $w\neq0$.  If $I_0$ were also
nonempty, then for every $i\in I_0$, the $i$-th coordinate of
$Aw=\lambda_{\max}w$ would give
\[
 0=\sum_{j\in I_+}A_{ij}w_j.
\]
All terms in this sum are nonnegative, and every $w_j$, $j\in I_+$,
is strictly positive.  Hence
\[
 A_{ij}=0\qquad(i\in I_0,\ j\in I_+).
\]

By symmetry the entries in the opposite rectangular block vanish as
well.  Reordering the indices so that those in $I_0$ come first would
therefore put $A$ in the excluded block-diagonal form.  Thus
$I_0=\varnothing$, and every coordinate of $w$ is strictly positive.
\end{proof}

\begin{lemma}\label{lem:semidefinite}
The matrix $M$ is negative semidefinite; equivalently,
\[
 x^TMx\leq0\qquad\text{for every }x\in\R^m.
\]
\end{lemma}

\begin{proof}
Lemma~\ref{lem:seminegative-M} gives
\[
 a^TMa\leq0\qquad(a\in\mathbb Z_{\geq0}^m).
\]
Multiplying a nonnegative rational vector by a common denominator gives
the same inequality for $a\in\mathbb Q_{\geq0}^m$; continuity then
gives it for every $a\in\R_{\geq0}^m$.

Let $\lambda_{\max}$ be the largest eigenvalue of $M$.
Since $M_{ij}=\int_{D_i\cap D_j}\rho\geq 0$ when $i\neq j$, and Lemma~\ref{lem:sign-matrix}(i), there is a
nonzero vector $a\geq0$ such that
$Ma=\lambda_{\max}a$.  Therefore
\[
 \lambda_{\max}\lVert a\rVert^2=a^TMa\leq0.
\]
Thus $\lambda_{\max}\leq0$.  Since every other eigenvalue is at most
$\lambda_{\max}$, all eigenvalues of $M$ are nonpositive, which is
exactly negative semidefiniteness.
\end{proof}
The main result in this section is the following 
\begin{thm}\label{thm:negative}
The matrix $M$ is negative definite.
\end{thm}
\begin{proof}
By Lemma~\ref{lem:semidefinite}, the only possible failure of negative
definiteness is a zero eigenvalue.  Suppose such an eigenvalue exists. After a simultaneous permutation of its rows and columns, write $M$
as a block-diagonal matrix
\[
 M=\operatorname{diag}\bigl(M^{(1)},\ldots,M^{(q)}\bigr),
\]
where none of the blocks $M^{(\nu)}$ admits a further decomposition
into two nonempty diagonal blocks.   Since $M$ has eigenvalue $0$, at least one block has
eigenvalue $0$.  Fix such a block.  It is negative semidefinite:
indeed, a vector supported on this block can be extended by zero to a
vector $x\in\R^m$, and then its quadratic form is $x^TMx\leq0$ by
Lemma~\ref{lem:semidefinite}.  Hence the largest eigenvalue of this
block is $0$.
Lemma~\ref{lem:sign-matrix}(ii), applied to this block, gives numbers
$a_i>0$ such that
\[
 \sum_jM_{ij}a_j=0
\]
for every index in this block.  Extend $a$ by zero on the other blocks;
then $Ma=0$ on all indices.  Put
\[
 E:=\sum_i a_iD_i,
\]
where the sum is over the block.
The crucial consequence, used repeatedly below, is
\begin{equation}\label{eq:E-orthogonal}
 \rho\cdot E\cdot D_i=(Ma)_i=0
 \qquad\text{for every }D_i\text{ occurring in }E.
\end{equation}

\smallskip
\noindent\emph{Step 1: the asymptotic J-slope inequality gives
$E^3\geq0$.}
Lemma~\ref{lem:seminegative-M} applies to divisors with integer coefficients,
whereas the coefficients $a_i$ of $E$ are real.  We therefore
approximate $E$ at a large scale.  For each large integer $k$, set  $m_i^{(k)}:=\lfloor ka_i\rfloor$, and set
\[
 F_k:=\sum_i m_i^{(k)}D_i=kE+\Delta_k,
\]
where the coefficients of $\Delta_k$ are uniformly bounded.  Since
\eqref{eq:E-orthogonal} gives
\begin{equation}\label{eq:rounded-asymptotics}
\begin{aligned}
 \rho F_k^2
 &=k^2\rho E^2+2k\rho E\Delta_k+\rho\Delta_k^2
   =\rho\Delta_k^2=O(1),\\
 F_k^3&=k^3E^3+O(k^2).
\end{aligned}
\end{equation}
The first equality is due to the useful cancellation: both $\rho E^2$ and
$\rho E\Delta_k$ are linear combinations of the zero numbers in
\eqref{eq:E-orthogonal}.

We next justify using one scale $\lambda=s/k$ for all large $k$.
The construction in Lemma~\ref{lem:seminegative-M} needs exactly two
facts: the class $\alpha-\lambda F_k$ is K\"ahler, and the
reference asymptotic J-slope $\lim_{R\to\infty}R^{-1}\cJ_{\omega}(\varphi_R)$ is positive.
For the first fact, observe that
\[
 \alpha-\frac{s}{k}\{F_k\}
 =\alpha-s\{E\}-\frac{s}{k}\{\Delta_k\}.
\]
The classes $\{\Delta_k\}$ stay in a bounded subset of the
finite-dimensional space $H^{1,1}(X,\R)$.  Choose $s_0>0$ so small
that $\alpha-sE$ remains in a fixed compact subset of the K\"ahler
cone for $0\leq s\leq s_0$.  After decreasing $s_0$, the displayed
class is K\"ahler for every large $k$ and every $0<s<s_0$.

For the second fact, set $G_k:=F_k/k$, so $G_k\to E$ as $k\to \infty$.
The exact asymptotic J-slope formula
\eqref{eq:asymptotic-omega-slope-polynomial} gives
\begin{equation}\label{eq:uniform-reference-weight}
\begin{aligned}
 k\lim_{R\to \infty}\frac{\cJ_\omega(\varphi_R)}{R}
  =\frac14\alpha^2G_k\,s^2
       -\frac13\alpha G_k^2\,s^3
       +\frac18G_k^3\,s^4.
\end{aligned}
\end{equation}
Because $E$ is a nonzero effective combination of divisors,
$\alpha^2E>0$.  Since $G_k\to E$, the first coefficient on the
right of \eqref{eq:uniform-reference-weight} is bounded below by a
positive number, while the other
coefficients remain bounded.  We may therefore fix
$0<s<s_0$, sufficiently small and independent of $k$, so that the
right-hand side of \eqref{eq:uniform-reference-weight} is positive for
every large $k$.
Thus $\lambda=s/k$ lies in the full admissible range used in
Lemma~\ref{lem:seminegative-M}: the logarithmic potentials are K\"ahler and their
reference asymptotic $J$-slope is positive.

Use \eqref{eq:seminegative} for $F_k$ at $\lambda=s/k$.  By
\eqref{eq:rounded-asymptotics},
\[
\begin{aligned}
 0
 &\leq-\rho F_k^2\frac{s^3}{k^3}
       +\frac{c_{\alpha,\beta}}{4}F_k^3\frac{s^4}{k^4}\\
 &=O(k^{-3})+\frac{c_{\alpha,\beta}s^4}{4k}E^3+O(k^{-2}).
\end{aligned}
\]
Multiplying by $k$ and letting $k\to\infty$ gives
$0\leq(c_{\alpha,\beta}s^4/4)E^3$.  Since $c>0$, this proves $E^3\geq0$.

\smallskip
\noindent\emph{Step 2: the surface Hodge index theorem gives
$E^3\leq0$.}
For every $D_i$ occurring in $E$, choose an embedded resolution
$\nu_i:\widehat D_i\to D_i$, and set
\[
 L_i:=\nu_i^*(\rho|_{D_i}),\qquad
 N_i:=\nu_i^*(E|_{D_i}).
\]
The embedded resolution is obtained by blowing up smooth centers in
$X$; hence $\widehat D_i$ is a smooth compact K\"ahler surface.
Because $\rho$ is nef, $L_i$ is nef.  It has positive square because
$D_i$ is a null divisor:
\[
 L_i^2=\rho^2D_i=\beta^2D_i>0
\]
(and $\beta^2D_i$ is the positive $\beta$-volume of $D_i$).
Moreover, by \eqref{eq:E-orthogonal},
\[
 L_iN_i=\rho ED_i=(Ma)_i=0.
\]
By the Hodge index theorem on $\widehat{D}_i$, the intersection form has signature $(1,h^{1,1}-1)$. Since $L_i^2>0$, it's restricion to $L_i^\perp$ is negative definite. 
Equality is possible only for  $N_i=0$ in $H^{1,1}(\widehat D_i,\R)$.
Since $E=\sum_i a_iD_i$, we may compute the cubic intersection:
\[
 E^3=\sum_i a_iE^2D_i=\sum_i a_iN_i^2\leq0.
\]
Step~1 gave the reverse inequality.  Since every $a_i>0$, equality of
the sum forces 
\begin{equation}\label{eq:restrictionzero}
 E^3=0,\qquad N_i=0\quad\text{for every }i.
\end{equation}
In ordinary words, the cohomology class of $E$ restricts to zero on
the resolution of each divisor occurring in $E$.

\smallskip
\noindent\emph{Step 3: equality would force $E$ to be nef.}
Fix any K\"ahler class $h$.  We shall check, dimension by dimension,
that for every positive-dimensional irreducible analytic subvariety
$V\subset X$, and every $s>0$,
\begin{equation}\label{eq:DPpositive}
 (E+sh)^{\dim V}\cdot V>0.
\end{equation}

\emph{Curves.}
Let $C$ be an irreducible curve.  If
$C\subset\operatorname{Supp}E$, choose $D_i\supset C$.  There is an
irreducible curve $\widehat C\subset\widehat D_i$ mapping generically
finitely onto $C$.  
\eqref{eq:restrictionzero} give
\[
 \deg(\widehat C/C)\,E\cdot C=N_i\cdot\widehat C=0.
\]
If $C\not\subset\operatorname{Supp}E$, effectiveness of $E$ gives
$E\cdot C\geq0$.  In both cases,
\[
 (E+sh)\cdot C=E\cdot C+s\,h\cdot C>0.
\]

\emph{Surfaces.}
Let $S$ be an irreducible surface.  If
$S\subset\operatorname{Supp}E$, then $S=D_i$ for some $i$.
The equality $N_i=0$ and the projection formula give
$EhD_i=E^2D_i=0$, and hence
\[
 (E+sh)^2D_i=s^2h^2D_i>0.
\]
If $S\not\subset\operatorname{Supp}E$, then $E\cdot S$ is an
effective curve cycle supported on $\operatorname{Supp}E$.  The curve
case just proved says that $E$ has zero intersection with every
component of this cycle.  Therefore
$E^2S=0$, while $EhS\geq0$, and
\[
 (E+sh)^2S=s^2h^2S+2s\,EhS>0.
\]

\emph{The threefold $X$.}
By \eqref{eq:restrictionzero},
\[
 E^2h=\sum_i a_i\,EhD_i=0,
\]
and Step~2 gives $E^3=0$.  Consequently
\[
 (E+sh)^3=s^3h^3+3s^2Eh^2>0.
\]
This proves \eqref{eq:DPpositive} in every possible dimension.

Let $\mathcal P$ denote the Demailly--P\u{a}un numerical positive
cone, consisting of classes $\gamma$ for which
$\gamma^{\dim V}\cdot V>0$ on every positive-dimensional irreducible
analytic subvariety $V$.  The preceding calculation says that the
entire connected ray
\[
 \{\,E+sh:s>0\,\}
\]
lies in $\mathcal P$.  For $s\gg1$, the class $E+sh$ is K\"ahler,
because $E+sh=s(h+s^{-1}E)$ and $h+s^{-1}E$ is a small perturbation
of $h$.  The Demailly--P\u{a}un component theorem says that the
K\"ahler cone is precisely one connected component of $\mathcal P$.
The whole ray must therefore lie in that component.  Thus $E+sh$ is
K\"ahler for every $s>0$, and its limit $E$ as $s\downarrow0$ is
nef.

Finally, $E$ is a nonzero class: indeed
$h^2E=\sum_i a_i h^2D_i>0$.  Every $D_i$ is null, so
\[
 z\cdot E=\sum_i a_i z\cdot D_i=0.
\]
But every nef class is modified nef, and Lemma~\ref{lem:MN} says that
$z$ is strictly positive on every nonzero modified-nef class.  This is
the required contradiction.  Hence $M$ is negative definite.
\end{proof}

\begin{cor}\label{cor:antinef}
There are positive rational numbers $a_i$ such that, for
\[
 E=\sum_{i=1}^m a_iD_i,
\]
one has
\begin{equation}\label{eq:antinef}
 \rho ED_j<0\qquad(1\leq j\leq m).
\end{equation}
    
\end{cor}
\begin{proof}
    Since $M$ is negative definite, it is invertible. The off-diagonal entries $M_{ij}=\int_{D_i}\rho\wedge\theta_{D_j}\geq 0$. Let $\mathfrak{1}=(1,1,...,1)^T$, and solve $Ma=-\mathfrak{1}, a\in \mathbb{R}^m$. We claim that every coordinate of $a$ is positive.  Write $a=a_+-a_-$, where $(a_+)_i=\max(a_i,0)$ and $(a_-)_i=\max(-a_i,0)$. If $a_-\neq 0$, since the support of $a_+$ and $a_-$ are disjoint and $M_{ij}\geq 0$ when $i\neq j$, then, $\langle Ma_+,a_-\rangle\geq 0$. But the negative definiteness gives $\langle Ma_-,a_-\rangle<0$. Hence \[
    \langle Ma,a_-\rangle=\langle Ma_+,a_-\rangle-\langle Ma_-,a_-\rangle>0.
    \]
    Since $Ma=-\mathfrak 1$, $\langle Ma,a_-\rangle=-\sum_{i=1}^m(a_-)_i <0$, a contradiction. Thus $a_i\geq 0.$ Moreover, if $a_i=0$ for some $i$, then the $i$-th equation of $Ma=-\mathfrak 1$ gives \[
    -1=(Ma)_i=\sum_{j\neq i}M_{ij}a_j\geq0.
    \]We get another contradiction. Therefore $a_i>0$ for all $i=1,...,m$.

    We choose $E=\sum_ia_iD_i$ for $a$ being the solution to $Ma=-\mathfrak 1$. Then $\rho ED_j=(Ma)_j=-1<0$. Positivity of coordinates and all these inequalities are open conditions. Approximating $a$ by rational vectors proves the corollary. 
\end{proof}

\section{Strict positive perturbation on a resolution
}
Put
\[
 C:=C_1\cup\cdots\cup C_\ell=\operatorname{Null}(\rho).
\] 
\begin{lemma}\label{lem:resolved-model}
There exists a resolution
\[
 \mu:Y\longrightarrow X,
\]
which is an isomorphism over $X\setminus C$, effective
$\mu$-exceptional integral divisors $H,Q,F$, supported over $C$, positive integers $p,q$, a K\"ahler current $T\in \rho$, a smooth closed form $\theta$, a
K\"ahler form $h_Y$ on $Y$, and  numbers $\gamma,\delta,b,u_0>0$ such that  \[
\text{Supp }H=\text{ Supp }F=\mu^{-1}(C)\] and 
\begin{enumerate}[label=\textup{(\roman*)}]
\item $T$ has analytic singularities exactly along $C$, and $\mu^*T=\theta+\gamma [H], $ with $\theta\geq \delta\mu^*h$.
\item $F=qH+pQ$.
\item For every $0<u<u_0$,  the class $\mu^*\rho-\mu[F]$ containsa smooth K\"ahler form $\kappa_u$ with $\kappa_u\geq buh_Y$.
\end{enumerate}
\end{lemma}
 
\begin{proof}
If $C=\emptyset$, then $E_{\rm nK}(\rho)=\emptyset$, so
$\rho$ is K\"ahler by the theorem of Collins-Tosatti\cite{CT15}. Take \[
\mu=\text{id}_X, \qquad H=Q=F=0, p=q=1,
\]choose K\"ahler form $T\in \rho$, put $\theta=T,\gamma=1,h_Y=h$ and choose $\delta>0$ with $T>2\delta h$. Set $\kappa_u:=T$. The lemma follows.  Hence, we assume that $C\neq \emptyset.$ 
By lemma~\ref{lem:rho-nef}, the class $\rho$ is nef and
$\rho^3>0$.  Collins--Tosatti's Theorems~1.1 and~2.2 give a
K\"ahler current $
 T\in\rho$, with analytic singularities and $ E_+(T)=C$. Choose it so that $T>2\delta h$ for some $\delta>0$, and let $\cI_T\subset \cO_X$ be the ideal sheaf defining the analytic singularities of $T$. Locally the potential of $T$ can be written as \[
 \gamma\log\left(\sum_\ell|f_\ell|^2\right)+g, \qquad \gamma>0, g\in C^\infty,
 \]where the $f_\ell$ generate $\cI_T$.
 
 Take an embedded principalization of $\cI_T$ which is isomorphism outside $C$: $\mu:Y\to X$, $\cI_T\cdot\cO_Y=\cO_Y(-H)$.  Then 
\begin{equation}\label{eq:numerical-resolved-current}
 \mu^*T=\theta+\gamma[H],
 \qquad
 \theta\geq\delta\,\mu^*h.
\end{equation}
Here $\theta$ is smooth and closed, and
$H$ is an effective integral $\mu$-exceptional divisor with $\supp  H=\mu^{-1}(C)$. To verify the lower bound of $\theta$, the positive current $T-\delta h$ has the same divisorial Lelong number as $T$. Its pullback therefore has Siu's decomposition $\mu^*(T-\delta h)=\theta-\delta\mu^*h+\gamma [H]$. The current $\theta-\delta\mu^*h$ is positive by Siu's decomposition theorem, which proves the desired inequality for $\theta$. 
We also use the standard exceptional correction for a composition of
blow-ups: there are an effective integral $\mu$-exceptional divisor
$Q$ with $\supp Q\subset \mu^{-1}(C)$ and a number $t_0>0$ such that
\begin{equation}\label{eq:exceptional-kahler-class}
 \mu^*\{h\}-t{Q}
 \quad\text{is a K\"ahler class for }0<t\leq t_0.
\end{equation}
Fix a K\"ahler form $h_0=\mu^*h-t_0\theta_{Q}\in \mu^*h-t_0[Q]$, then for any $t\in (0,t_0)$, \[\mu^*h-t\theta_Q=(1-\frac{t}{t_0})\mu^*h+\frac{t}{t_0}h_0\geq \frac{t}{t_0}h_0.\] 
Choose positive integers $p,q$ such that $0<\frac pq<\frac{\delta t_0}{\gamma}$, and put $
 \lambda:=\gamma\frac pq$, $ F:=qH+pQ$, and $\gamma_*:=\frac{\gamma}{q}$.
Then $0<\lambda/\delta<t_0$, so
\begin{equation}
 \kappa
 :=\theta-\lambda\theta_Q \notag=(\theta-\delta\mu^*h)
   +\delta\left(
      \mu^*h-\frac{\lambda}{\delta}\theta_Q
    \right)>0.       \label{eq:fixed-kappa-numerical}
\end{equation}  Thus $\kappa$ is a smooth
K\"ahler form.  Moreover,
\[
 \gamma H+\lambda Q
 =\frac{\gamma}{q}(qH+pQ)
 =\gamma_*F,
\]
and
\begin{equation}\label{eq:fixed-kappa-class}
 \{\kappa\}=\mu^*\rho-\gamma_*\{F\}.
\end{equation}
The divisor $F$ is effective, integral, and $\mu$-exceptional.  It
is supported over $C$, because both $H$ and $Q$ are supported
there.  Since $q>0$ and $\operatorname{Supp}H=\mu^{-1}(C)$, it also
satisfies $\operatorname{Supp}F=\mu^{-1}(C)$.

It remains to replace the fixed coefficient $\gamma_*$ by an
arbitrary sufficiently small coefficient $u$.  Since $\rho$ is
nef, for every $s>0$ there is a smooth closed form
\[
 \alpha_s\in\rho,\qquad \alpha_s\geq-sh.
\]
Fix a K\"ahler form $h_Y$ on $Y$.  By compactness, choose constants
$b_0,A>0$ such that
\begin{equation}\label{eq:numerical-fixed-comparisons}
 \kappa\geq b_0h_Y,\qquad \mu^*h\leq Ah_Y.
\end{equation}
For $0<u<\gamma_*$, define
\begin{equation}\label{eq:numerical-kappa-u}
 \kappa_u
 :=
 \frac{u}{\gamma_*}\kappa+
 \left(1-\frac{u}{\gamma_*}\right)\mu^*\alpha_{u^2}.
\end{equation}
Using \eqref{eq:fixed-kappa-class}, we calculate
\[
\begin{aligned}
 \{\kappa_u\}
 &=
 \frac{u}{\gamma_*}
   \bigl(\mu^*\rho-\gamma_*\{F\}\bigr)
 +\left(1-\frac{u}{\gamma_*}\right)\mu^*\rho\\
 &=\mu^*\rho-u\{F\}.
\end{aligned}
\]
Furthermore, $\mu^*\alpha_{u^2}\geq-u^2\mu^*h$, and hence
\eqref{eq:numerical-fixed-comparisons} gives
\[
\begin{aligned}
 \kappa_u
 &\geq
 \frac{b_0}{\gamma_*}u\,h_Y
 -\left(1-\frac{u}{\gamma_*}\right)u^2\mu^*h\\
 &\geq
 \left(\frac{b_0}{\gamma_*}u-Au^2\right)h_Y.
\end{aligned}
\]
After decreasing $u_0>0$, we may arrange that
\[
 0<u_0<\gamma_*,
 \qquad
 Au_0\leq\frac{b_0}{2\gamma_*}.
\]
It follows that, for $0<u<u_0$,
\[
 \kappa_u\geq
 \frac{b_0}{2\gamma_*}u\,h_Y.
\]
Taking $b=b_0/(2\gamma_*)$ proves the lemma. 
\end{proof}

\begin{rmk}[Why two exceptional divisors are retained]
\label{rem:two-resolved-divisors}
The analytic divisor $H$ and the numerical correction divisor $F$
play different roles.  The resolution of the fixed current $T$
gives the exact current identity
$\mu^*T=\theta+\gamma[H]$, whereas the blow-up correction gives only
the K\"ahler classes $\mu^*\rho-u\{F\}$, with
$F=qH+pQ$.   In the descent argument
below, $H$ controls the analytic singularities and $F$ controls the
numerical perturbation; the relation $F=qH+pQ$ gives the strict
comparison needed to glue the two currents.
\end{rmk}

We next determine the null divisors upstairs on $Y$ similar to lemma \ref{lem:finite-null-divisor}. Let $\overline{\Eff}(Y)$ denote the closed cone generated by effective $\R$-divisors on $Y$.

\begin{lemma}\label{lem:finite-generated-null-blowup}Let $\cG$ be the finite collection consisting of all $\mu$-exceptional prime divisors and the strict transform of $\tilde{D_i}$ for J-null prime divisors $D_i$ on $X$. Then \[
\{\xi\in \xi\in\overline{\Eff}(Y):\mu^*z\cdot\xi=0\}=\sum_{G\in \cG}\R_{\ge 0}\{G\}
\]
\end{lemma}
\begin{proof}  It's clear that any $\R$-combinations of $\mu$-exceptional divisors and strict proper transforms of $J$-null divisors lie in the zero set of $\mu^*z$. 
Conversely, let $\xi\in \overline{\Eff}(Y)$ satisfies $\mu^*z\cdot \xi=0$. Choose effective $\R$-divisors whose class converge to $\xi$, and pass to a weak limit of their integration currents. Their mass against $h_Y^{n-1}=h_Y^2$ are bounded uniformly by cohomology, so it gives a current $T\in \xi$. Then we have \[
z\cdot \{\mu_*T\}=\mu^*z\cdot \xi=0.
\]
By lemma \ref{lem:finite-null-divisor}, \[
\{\mu_*T\}=\sum_ib_i\{D_i\}, \qquad b_i\geq 0.
\]
The null divisors form an exceptional family, so Boucksom rigidity theorem \cite[Proposition 3.13]{B04} implies this class has the unique positive representative \[
\mu_*T=\sum_ib_i[D_i].
\]
On $Y\setminus \mu^{-1}(C\cup_iD_i)$, the map $\mu$ is a biholomorphism and the preceding identity forces $T=0$ there. The support theorem for positive closed $(1,1)$-currents therefore gives \[
T=\sum_{G\in \cG}b_G[G], \qquad b_G\geq 0.
\]

\end{proof}

We now show the existence of direction of strict positive deformation of $J$-slope. 
\begin{lemma}\label{lem:strict-cone-blowup}Let $E=\sum_ia_iD_i$ be the divisor in corollary~\ref{cor:antinef}.
There is a positive rational constant $A_0$ and there are sufficiently
small rational numbers $s,\varepsilon,\eta >0$,
such that, for
\[
 \widehat r
 :=\mu^*\rho-s\{F\}-A_0s^2\mu^*\{E\}-\cab\varepsilon\mu^*\beta, \qquad \widehat\beta:=\mu^*\beta+\eta \{h_Y\},
\]
one has:
\begin{enumerate}[label=\textup{(\roman*)}]
\item $\widehat r$ and $\widehat\beta$ are K\"ahler classes;
\item $(\widehat r^2-\widehat\beta^2)\cdot \xi>0, \forall 0\neq  \xi\in\overline{\operatorname{Eff}}(Y);$
\item $\widehat r^3-3\widehat r\,\widehat\beta^2>0.$
\end{enumerate}
\end{lemma}
\begin{proof}
    First omit the $\varepsilon$- and $\eta$-terms.  Put
\[
 \widehat r_s^0:=\mu^*\rho-s\{F\}-A_0s^2\{\mu^*E\}.
\]We suppress braces around divisor classes in the intersection computations below. 
For every irreducible divisor $D\subset Y$, $\mu^*z\cdot D=z\cdot\mu_*D\geq0.$
Moreover, we claim that 
\begin{equation}\label{eq:linear-exceptional-zero}
 \mu^*\rho\cdot F\cdot D=0.
\end{equation}
Indeed, if $D$ is not contained in any component of $F$, then $D\cdot F$ is an effective circle supported over $C$. It follows that $\mu^*\rho\cdot D\cdot F=0$ since $\rho\cdot C=0$. If $D$ is one the irreducible component of $F$. We define a current $T_D:=\mu_*(\theta_D\wedge[D])$, where $\theta_D$ if the curvature form of $\cO(D)$ with a hermitian metric. The current $T_D$ is closed, normal(of order $0$), of bidegree $(2,2)$, and it's supported on $\mu(D)$. If $\mu(D)$ is a point, then $T_D=0$. Otherwise $\mu(D)=C_i$ for some null curve $C_i$ of $\rho$.  By support theorem \cite[Corollary 2.14]{D12}, $T_D=c[C_i]$. Therefore \[
\mu^*\rho \cdot D^2=\int_Y\mu^*\rho\theta_D^2=\int_Y\mu^*\rho \theta_D\wedge[D]=\int_X\rho \wedge T_D=c\int_{C_{i,\reg}}\rho=0.
\]
It follows that $\mu^*\rho \cdot D\cdot F=0.$

Consequently, for every
$\xi\in\overline{\operatorname{Eff}}(Y)$,
\begin{equation}\label{eq:second-expansion}
\begin{aligned}
 &\bigl((\widehat r_s^0)^2-(\mu^*\beta)^2\bigr)\cdot\xi=\mu^*z\cdot\xi+s^2L\cdot\xi+s^3R_s\cdot\xi,
\end{aligned}
\end{equation}
where
\[
 L:=F^2-2A_0\mu^*\rho\,\mu^*E,
 \qquad
 R_s:=2A_0F\cdot\mu^*E+A_0^2s(\mu^*E)^2.
\]
We show that $L$ is strictly positive on $\{\xi\in\overline{\operatorname{Eff}}(Y):\mu^*z\cdot\xi=0\}$.

If $D$ is a $\mu$-exceptional irreducible divisor, then
Lemma~\ref{lem:resolved-model} makes
$\mu^*\rho-sF$ K\"ahler.  Since $(\mu^*\rho)^2D=\mu^*\rho F\cdot D=0,$
we obtain
\[
 0<(\mu^*\rho-sF)^2\cdot D=s^2F^2\cdot D,
\]
and hence $F^2\cdot D>0$.  Also
\[
 \mu^*\rho\cdot\mu^*E\cdot D=\int_Y\mu^*\rho \mu^*\theta_E\wedge[D]=\int_X\rho\theta_E\mu_*[D]=0,
\]
where $\theta_E$ is the curvature of $\cO(E)$ and the last equality holds due to similar argument using the support theorem and the fact $\mu(D)$ has codimension at least $2$.
Thus $L\cdot D=F^2\cdot D>0 $ on every exceptional generator.

If $D=\widetilde D_i$ is the strict
transform of a null divisor, then
\[
 L\cdot D=F^2\widetilde D_i-2A_0\rho ED_i.
\]
Here we used that \[
2A_0\mu^*\rho \cdot \mu^*E\cdot \widetilde{D}_i=2A_0\int_{\widetilde{D}_{i,\reg}}\mu^*\rho\, \mu^*\theta_E=2A_0\int_{D_{i,\reg}}\rho\theta_E=2A_0\rho ED_i.
\]
The last intersection is negative by
Corollary~\ref{cor:antinef}; hence one fixed sufficiently
large rational $A_0$ makes $L\widetilde D_i>0$ for every $i$.

We now give the uniform compactness argument explicitly.  Normalize the
closed effective-divisor cone by
\[
 B:=\{\xi\in\overline{\operatorname{Eff}}(Y):
                  h_Y^2\cdot\xi=1\}.
\]
This set is compact.  Indeed, it is closed.  If it were unbounded for
some norm on the finite-dimensional space $H^{1,1}(Y,\mathbb R)$,
we could choose $\xi_j\in B$ with $\|\xi_j\|\to\infty$.  After
normalizing and taking a subsequence,
\[
 \frac{\xi_j}{\|\xi_j\|}\longrightarrow\xi_\infty
 \in\overline{\operatorname{Eff}}(Y),\qquad
 \|\xi_\infty\|=1,
\]
while $h_Y^2\cdot\xi_\infty=0$.  This is impossible: a nonzero
pseudo-effective class has strictly positive mass against
$h_Y^2$.  Concretely, choose effective divisor currents whose classes
converge to $\xi_\infty$.  Their $h_Y^2$-masses are bounded, so a
subsequence converges weakly to a positive current $T_\infty$ in
$\xi_\infty$.  Since $\|\xi_\infty\|=1$, this current is nonzero,
and hence
$\int_YT_\infty\wedge h_Y^2>0$, a contradiction.

The functional $\mu^*z$ is nonnegative on all of $B$ by $J$-nefness.  Let
\[
 B_0:=\{\xi\in B:\mu^*z\cdot\xi=0\}.
\]
If $B_0=\varnothing$, then $\mu^*z$ has a positive minimum on $B$,
and \eqref{eq:second-expansion} is positive for all sufficiently small
$s$.  Suppose $B_0\neq\varnothing$.  Lemma \ref{lem:finite-generated-null-blowup} and the strict positivity of $L$ on every generator show that
\[
 a:=\min_{\xi\in B_0}L\cdot\xi>0.
\]
Choose a neighborhood $U$ of $B_0$ in $B$ such that
\[
 L\cdot\xi\geq\frac a2\qquad(\xi\in U).
\]
The compact set $B\setminus U$ contains no zero of $\mu^*z$, and
$\mu^*z\geq0$ everywhere; hence
\[
 b:=\min_{\xi\in B\setminus U}\mu^*z\cdot\xi>0.
\]

The family $R_s$ is bounded for $0\leq s\leq1$.  Since $B$ is
compact and the intersection pairings are continuous, there is
$M>0$ such that, for $0<s\leq1$ and $\xi\in B$,
\[
 |L\cdot\xi|+|R_s\cdot\xi|\leq M.
\]
For $\xi\in U$, we obtain
\[
 (\mu^*z+s^2L+s^3R_s)\cdot\xi
 \geq s^2\left(\frac a2-Ms\right)>0
\]
when $s$ is small.  For $\xi\in B\setminus U$,
\[
 (\mu^*z+s^2L+s^3R_s)\cdot\xi
 \geq b-Ms^2-Ms^3>0.
\]
Thus
\[
 \bigl((\widehat r_s^0)^2-(\mu^*\beta)^2\bigr)\cdot\xi>0
 \quad
 (0\neq\xi\in\overline{\operatorname{Eff}}(Y))
\]
for all sufficiently small $s>0$.

By Lemma~\ref{lem:resolved-model},
$\mu^*\rho-sF$ has a K\"ahler representative bounded below by a
constant times $s h_Y$.  The additional term
$-A_0s^2\mu^*E$ is of smaller order, so
$\widehat r_s^0$ is K\"ahler for small $s$.
The top-degree inequality follows from $\rho^3-3\rho\beta^2=2\beta^3>0$ by
continuity.  Finally choose $0<\varepsilon\ll s^2$, and then
$0<\eta\ll s^2$.  The square-cone margin just obtained is of order
$s^2$, so these are the precise scales needed for that inequality.  The top inequality is open at $s=0$.  Thus all three
strict properties persist, proving (i)--(iii).
\end{proof}

\begin{lemma}[Fang-Ma]
    Let $r,\beta$ be K\"ahler classes on a connected compact K\"ahler
threefold, and let $\chi\in\beta$ be K\"ahler.  Assume
\[
 (r^2-\beta^2)\cdot D>0
 \quad\text{for every prime divisor }D,
 \qquad
 r^3-3r\beta^2>0.
\]
There is a K\"ahler form $R\in r$ such that
\[
 R^2-\chi^2>0.
\]
\end{lemma}
\begin{proof}
    Put $f:=\frac{r^3-3r\beta^2}{\beta^3}>0$, and use Fang--Ma's notation exactly as follows.  For a K\"ahler form
$R$, they put
\[
 \Omega(R)=\exp R=1+R+\frac{R^2}{2}+\frac{R^3}{6}
\]
and study
\[
 \kappa\,\Omega(R)^{[3]}=(\Lambda\wedge\Omega(R))^{[3]}
 \tag{FM 1.1}
\]
\cite[equation~(1.1)]{FM24}.  Take $\kappa=1$ and
\[
 \Lambda^{[1]}=0,\qquad
 \Lambda^{[2]}=\frac{\chi^2}{2},\qquad
 \Lambda^{[3]}=\frac{f\chi^3}{6}.
\]
Let $\mathring\Lambda=\Lambda^{[1]}+\Lambda^{[2]}$.  Then
$\mathring\Lambda=\chi^2/2$ is (2)-uniformly positive in the
sense of Fang--Ma, with reference form $\chi$ and uniform constant
$m=1$.  The top-degree component
$\Lambda^{[3]}=f\chi^3/6$ is strictly positive, and hence is an
almost-positive top form for the same data.  Thus their structural
hypothesis  \textbf{H1} is satisfied
\cite{FM24}.
Then (FM 1.1) is, with no suppressed factorials,
\begin{equation}\label{eq:FM-specialized}
 \frac{R^3}{6}
 =\frac{\chi^2\wedge R}{2}+\frac{f\chi^3}{6},
 \qquad\text{that is}\qquad
 R^3=3\chi^2\wedge R+f\chi^3.
\end{equation}
Its cohomological equality is
\[
 \int_X\frac{r^3}{6}
 =\int_X\frac{\beta^2r}{2}+\frac f6\int_X\beta^3.
\]
Fang--Ma's cone condition is their equation~(1.5):
\[
 \bigl((1-\Lambda)\wedge\exp R\bigr)^{[2]}
 =\frac12(R^2-\chi^2)>0.
 \tag{FM 1.5}
\]
Their Definition~1.6 says that $r$ is
$([\Lambda],1)$-positive when the degree-$\dim Y$ part of
$[\exp r]\,[1-\Lambda]$ has positive integral on every proper
irreducible subvariety $Y$, and has nonnegative integral on $X$.
In the present specialization these conditions are exactly
\[
 r\cdot C>0
 \quad\text{for every curve }C,
\]
and
\[
 \frac12(r^2-\beta^2)\cdot D>0
 \quad\text{for every irreducible divisor }D.
\]
The first follows because $r$ is K\"ahler, the second is one of the hypotheses of the lemma, and the condition on $X$ is equality by
the displayed cohomological identity. Fang-Ma's theorem 1.10 \cite{FM24} gives a K\"ahler form $R\in r$ solving \eqref{eq:FM-specialized}. We may check that $R^2-\chi^2>0$ from the equation: If $\lambda_1,\lambda_2,\lambda_3>0$ are the eigenvalues
of $\chi$ relative to $R$, then \eqref{eq:FM-specialized} says
\[
 1=\lambda_1\lambda_2+\lambda_1\lambda_3+
   \lambda_2\lambda_3+f\lambda_1\lambda_2\lambda_3.
\]
Since $f>0$, every $\lambda_i\lambda_j<1$, which is precisely
$R^2-\chi^2>0$.
\end{proof}

Combinining with lemma \ref{lem:strict-cone-blowup}, we get 
\begin{cor}\label{cor:strict-cone-blowup}Let \[ \widehat{\chi}:=\mu^*\chi+\eta h_Y\in \widehat{\beta}.\] There is a K\"ahler form $\widehat{R}\in \widehat{r}$ such that \[
\widehat{R}^2-\widehat{\chi}^2>0.
\]
\end{cor}

\section{Construction of the analytic cone current}
We now construct an element in $\cK_{J}^{\rm an}(\alpha,\beta)$. We begin with the construction of a current satisfying the cone condition with analytic singularity along the null curve $C$ first. We take $\widehat{R}$ from corollary \ref{cor:strict-cone-blowup}.

\begin{prop}
\label{prop:analytic-descent}
Let $ r_X:=\rho-A_0s^2\{E\}-\cab\varepsilon\beta$.
There is a closed positive current $R\in r_X$, with analytic
singularities exactly along $C$, such that $ R^2-\chi^2\geq0$ 
in the following local-convolution sense: on every coordinate ball,
for every constant-coefficient K\"ahler form
$\chi_0$ satisfying $\chi\geq\chi_0$, the ordinary convolution
$R_\delta$ on the smaller ball satisfies
\[
 R_\delta^2-\chi_0^2\geq0.
\]
\end{prop}

We observe that the cone condition in the sense of local convolution indeed allows a negligible singular set.
\begin{lemma}\label{lem:convolution-negligible-set}Let $O\subset \C^3$ be a open set, and $Z\subset O$ have Lebesgue measure zero. Let $R$ be a positive $(1,1)$-current whcih is smooth on $O\setminus Z$, and satisfying $R^2-\chi^2>0$ on $O\setminus Z$, where $\chi$ is a smooth K\"ahler form on $O$. Then $R^2-\chi^2\geq 0$ in the sense of local convolution sense.   
\end{lemma}
\begin{proof}
    The coefficients of $R$ are measures, and decompose it into the absolute continuous part and singular part, say $R=A+R^s$ with $A$ being the absolute continuous part. Then $A$ is a smooth off $Z$ and satisfies $A^2-\chi^2>0$ on $O\setminus Z$. We say that a $3\times 3$ matrix $M$ with entries $m_{i\bar j}$ being Radon measure is positive is for any $\xi=(\xi^i,\xi^2,\xi^3)\in \C^3$, the measure $\xi^i\overline{\xi^j}m_{i\bar j}$ is a positive Radom measure. It clear that the decomposition of measure $\xi^i\overline{\xi^j}R_{i\bar j}$ is $\xi^i\overline{\xi^j}R_{i\bar j}=\xi^i\overline{\xi^j}A_{i\bar j}+\xi^i\overline{\xi^j}R^s_{i\bar j}$. Hence both $A$ and $R^s$ are positive.  The local convolution of the singular part is \[
    R^s_\delta(x)=\sum_{i,j}\sqrt{-1}\left(\int\rho_\delta(x-y)d R^s_{i\bar j}(y)\right)dz_i\wedge d\bar{z}_j.
    \]
    Hence $\xi^i\overline{\xi^j}R^s_{\delta,i\bar j}\geq 0$ for any $\xi\in \C^3$. From the monotonicity of $P_{\chi_0}(\cdot)$, it remains to show $A_\delta$ satisfies $A_{\delta}^2-\chi_0^2\geq 0$. This follows from the convexity of $P_{\chi_0}(\cdot)$ and $A(x)^2\geq \chi_0^2$ almost everywhere on $O$. 
\end{proof}

Now we prove a similar statement in proposition~\ref{prop:analytic-descent} without prescribing analytic singularity. 
\begin{lemma}
\label{lem:raw-pushdown}
Let $r_X:=\rho-A_0s^2\{E\}\cab\varepsilon\beta$ and $R_0:=\mu_*\widehat R$.  Then $R_0\in r_X$ is a K\"ahler
current, it is smooth on $X\setminus C$, and
\[
 R_0^2-\chi^2>0,
\]
in the local-convolution sense.
\end{lemma}
\begin{proof}
    Since $\widehat{r}=r_X+s\{F\}$ and $F$ is $\mu$-exceptional, we have $[R_0]\in r_X$: indeed, for any smooth closed form $\varphi$ of degree $(n-1,n-1)$ on $X$, we have \[
    \int_X\varphi \mu_*\widehat{R}=\int_Y\mu^*\varphi \widehat{R}=\int_{Y}\mu^*\varphi( \mu^* r_X-s[F])=\int_X\varphi\,r_X-s\int_X\varphi \mu_*[F].
    \]
    $\mu_*[F]$ is supported on $\mu(F)$ but $\text{codim}_\C\,\mu(F)\geq 2$. Support theorem implies that $\mu_*[F]=0$. Then Poincare duality implies that $[R_0]\in r_X. $

    The current $R_0$ is smooth on $X\setminus C$. And $\mu^* R_0=\widehat{R}$ on $Y\setminus \mu^{-1}(C)$. Moreover, $\widehat{\chi}=\mu^*\chi+\eta h_Y\geq \mu^*\chi,$ hence $
    \widehat{R}^2-(\mu^*\chi)^2\geq 0$ on $Y\setminus \mu^{-1}(C)$. It implies $R_0^2-\chi^2>0$ on $X\setminus C$. Then lemma \ref{lem:convolution-negligible-set} implies $R_0^2-\chi^2\geq 0$ in the local convolution sense. 
\end{proof}

\begin{lemma}\label{lem:raw-pullback}
$\mu^*R_0=\hat{R}+s[F]$. 
\end{lemma}
\begin{proof}
  Recall that the pullback of the closed positive current $\mu^*R_0$ is defined locally by $\ddc (\varphi_{R_0}\circ\mu)$ where $\varphi_{R_0}$ is a local potential of $R_0$.   Since
$\widehat R\in \mu^*r_X-s\{F\}$, the positive current
$\widehat R+s[F]$ lies in $\mu^*r_X$ and agrees with
$\mu^*R_0$ away from $F$.  Let \[
U=:\mu^*R_0-(\widehat{R}+s[F]).\]
 Then $U$ is a closed current of order zero supported on the exceptional divisors and the class $\{U\}=0\in H^{1,1}(Y,\R)$. The support theorem therefore gives \[
 U=\sum_{i}b_i[E_i]
 \] where the $E_i$ are prime $\mu$-exceptional divisors.

 independence of exceptional
divisor classes forces every coefficient to vanish.  The required
independence holds because $\mu$ is a composition of blow-ups with
smooth centers: at each step the new exceptional class is the new
summand modulo the pullback of the preceding $H^{1,1}$. 
\end{proof}

\begin{lemma}\label{lem:smaller-model}
Put $\tau=qs/2$.  There is a K\"ahler current $G\in r_X$, with
analytic singularities exactly along $C$, such that
\[
 G\geq d\chi,
\] for some $d>0$, and on the resolution $Y$, the divisorial part
of $\mu^*G$ is exactly $\tau H$.
\end{lemma}
\begin{proof}

The integer $q$ is fixed once and for all in
Lemma~\ref{lem:resolved-model}, so we may take $s$ small enough that
\[
 0<\tau<\gamma,
 \qquad
 \tau\leq\frac{\delta}{2\gamma}.
\]
Since $\rho$ is nef, choose a smooth closed form
\[
 \alpha_{\tau^2}\in\rho,
 \qquad
 \alpha_{\tau^2}\geq-\tau^2h,
\]
and define
\begin{equation}\label{eq:scaled-fixed-current}
 T_\tau
 :=
 \frac{\tau}{\gamma}T+
 \left(1-\frac{\tau}{\gamma}\right)\alpha_{\tau^2}.
\end{equation}
First observe that $T\geq\delta h$.  
Consequently
\begin{equation}\label{eq:scaled-current-lower-bound}
 T_\tau
 \geq
 \left(\frac{\delta\tau}{\gamma}-\tau^2\right)h
 \geq\frac{\delta\tau}{2\gamma}h.
\end{equation}
Moreover $T_\tau\in\rho$, it has analytic singularities exactly
along $C$, and its pullback is
\begin{equation}\label{eq:scaled-current-pullback}
 \mu^*T_\tau
 =
 \left[
  \frac{\tau}{\gamma}\theta+
  \left(1-\frac{\tau}{\gamma}\right)\mu^*\alpha_{\tau^2}
 \right]
 +\tau[H].
\end{equation}
Equivalently, in local expression
\[
 \varphi_T
 =\gamma\log\!\left(\sum_\ell|f_\ell|^2\right)+g,
 \qquad g\in C^\infty,
\]
a potential of $T_\tau$ is
\[
 \tau\log\!\left(\sum_\ell|f_\ell|^2\right)+g_\tau,
 \qquad g_\tau\in C^\infty.
\]

Let $\theta_E$ be the curvature form of $\cO(E)$ equipped with a smooth hermitian metric, and set
\begin{equation}\label{eq:local-model-G}
 G:=T_\tau-A_0s^2\theta_E-\cab\varepsilon\chi.
\end{equation}
This current lies in $r_X$.  To check its strict positivity without
hiding any estimate, choose constants $M_E,M_\chi>0$ such that
\[
 -M_Eh\leq\theta_E\leq M_Eh,
 \qquad
 \chi\leq M_\chi h.
\]
Equations \eqref{eq:scaled-current-lower-bound} and
\eqref{eq:local-model-G} give
\[
 G\geq
 \left(
  \frac{\delta\tau}{2\gamma}
  -A_0M_Es^2-cM_\chi\varepsilon
 \right)h.
\]
Because $\tau=qs/2$ is of order $s$, whereas
$s^2=o(s)$ and $\varepsilon\ll s^2$, the parameters in
Lemma~\ref{lem:strict-cone-blowup} may be chosen so small that
\[
 A_0M_Es^2\leq\frac{\delta\tau}{8\gamma},
 \qquad
 cM_\chi\varepsilon\leq\frac{\delta\tau}{8\gamma}.
\]
Thus $G\geq\delta\tau(4\gamma)^{-1}h$.  In particular, after
comparing the two fixed K\"ahler forms $h$ and $\chi$, there is a
number $d>0$ such that
\begin{equation}\label{eq:G-positive-lower-bound}
 G\geq d\chi.
\end{equation}
Subtracting smooth forms in \eqref{eq:local-model-G} does not change
the singularities.  Hence $G$ has analytic singularities
exactly along $C$, and \eqref{eq:scaled-current-pullback} shows that
its divisorial part on $Y$ is exactly $\tau H$.
\end{proof}

Now we note that around a arbitrary point, one can construct smooth function $g$ with almost positive Levi form and $g=L|z|^2$ near that point. 
\begin{lemma}
\label{lem:common-trace}
Let $S$ be a set of finite points on $X$. For each $a\in S$, choose a coordinate ball $B_a$ with the balls disjoint and put $q_a(z):=|z_1|^2+|z_2|^2+|z_3|^2.$ 
Given positive numbers $L_a$ and $\eta>0$, after shrinking
the $B_a$ there is a real function $g:X\to \R$
\begin{enumerate}[label=\textup{(\roman*)}]
\item $g_X=L_aq_a$ near $a$;
\item $\ddc g_X\geq-\eta\chi$.
\end{enumerate}
\end{lemma}

\begin{proof}
It is enough to construct the function in each of the pairwise
disjoint balls.  Fix $a\in S$, identify a smaller ball with
$\{|z|<R_a\}\subset\mathbb C^3$, and let $K>2$.  Choose a smooth radial cutoff $\vartheta_{a,K}$ which is $1$ on $\{|z|\leq R_ae^{-K}\}$, vanishes on $\{|z|=R_a\}$, satisfies $0\leq \vartheta_{a,K}\leq 1$( for example, smoothing of $K^{-1}\log(R_a/|z|)$). With respect to the Euclidean metric, it may be chosen so that \[
|d\vartheta_{a,K}|\leq \frac{A_a}{K|z|}, \qquad |\ddc \vartheta_{a,K}|\leq \frac{A_a}{K|z|^2}, 
\]for a uniform constant $A_a$. Since $q_a=|z|^2, |dq_a|=O(|z|)$, and $\ddc q_a\geq 0$, one can verify directly that \[\begin{aligned}
 \ddc\bigl(L_a\vartheta_{a,K}q_a\bigr)
 &=L_a\vartheta_{a,K}\ddc q_a
   +L_aq_a\ddc\vartheta_{a,K}\\
 &\quad+L_a d\vartheta_{a,K}\wedge d^cq_a
       +L_a dq_a\wedge d^c\vartheta_{a,K}\\
 &\geq-\frac{A'_aL_a}{K}\,\omega_{\rm Euc}
 \geq-\frac{A''_aL_a}{K}\,\chi.
\end{aligned}\]
Here the last comparison is uniform after the coordinate ball has
been fixed.  Choose $K=K_a$ so large that
$A''_aL_a/K_a\leq\eta$, and extend
$L_a\vartheta_{a,K_a}q_a$ by zero outside $B_a$.  It is smooth
because the cutoff vanishes on an outer collar, and it equals
$L_aq_a$ near $a$.  The supports for different $a$'s are
disjoint, so their sum is a globally smooth function satisfying the
same lower bound.  This sum is the required $g_X$
\end{proof}

We apply this lemma to $S=C_{\sing}$. 

\begin{lemma}
\label{lem:compatible-local-branches}
There are finitely many coordinate balls $W'_j\Subset W_j$ with $W'_j$  covering $C$, and smooth
real functions $\psi_j\in C^\infty(W_j)$ with the following
properties:
\begin{enumerate}[label=\textup{(\roman*)}]
\item on $W_j\setminus C$,
      \[
       B_j:=G+\ddc\psi_j>0,\qquad B_j^2-\chi^2>0;
      \]
\item on every overlap,
      \[
       \psi_i|_C=\psi_j|_C
       \quad\text{on }C\cap W_i\cap W_j;
      \]
\item there are numbers $a,\kappa>0$, independent of $j$, such that
      \begin{equation}\label{eq:local-branch-uniform-margins}
       B_j\geq a\chi,\qquad
       P_\chi(B_j+\chi)\leq1-2\kappa
       \quad\text{on }W_j\setminus C.
      \end{equation}
\end{enumerate}
\end{lemma}

\begin{proof}  Let $d>0$ be given by
\eqref{eq:G-positive-lower-bound}.  
Choose pairwise disjoint coordinate balls about the point $a\in C_{\sing}$, and $q_a=|z|^2$ there. There is $v_a>0$ such that $\ddc q_a\geq v_a\chi$. Choose $L_a$ so large that $d+L_av_a>1$ and apply lemma \ref{lem:common-trace} with $S=C_\sing$, $\eta=\frac{d}{4}$. We ostain $g_X\in C^\infty(X)$ such that \[
\ddc g_X\geq -\frac{d}{4}\chi, \qquad g_X=L_aq_a\quad \text{ near every }a\in C_\sing.
\]
Choose nested balls $W_a'\Subset W_a$ inside the regions where $g_a=L_aq_a$ holds. On $W_a$, set $\psi_a=g_X$. On $W_a\setminus C$, \[
G+\ddc\psi_a=G+L_a\ddc q_a\geq (d+L_av_a)\chi>\chi. 
\] So the desired strict inequalities hold on every singular-point chart.

The compact set $C\setminus \bigcup_{a\in C_\sing}W_a'$ is contained in $C_\reg$, and hence is covered by finitely many product coordinates balls \[
(z_1,z_2,w),\quad C=\{z_1=z_2=0\}.
\] 
Choose smaller product balls still covering this compact set. On the $j$-th ball, put $\Xi_j:=\ddc(|z_1|^2+|z_2|^2),$ and $\psi_j=g_X+N_j(|z_1|^2+|z_2|^2)$. The two nonzero eigenvalues of $\chi^{-1}\Xi_j$ are bounded below by a positive number $v_j>0$. We have $G+\ddc g_X\geq \frac{3}{4}d\chi$. The min-max principle therefore gives, for the ordered eigenvalues of $\chi^{-1}B_j$ satisfies \[
\lambda_1\geq \frac{3}{4}d,\qquad \lambda_2\geq \frac{3}{4}d+N_jv_j.
\]
We choose $N_j$ so that $\frac{3d}{4}\left(\frac{3d}{4}+N_jv_j\right)>1$. Hence $B_j^2-\chi^2>0$.

Note that on every ball $W_j$, normal quadratic term vanishes on $C$, while on every singular-point chart $W_a$, $\psi_a=g_X$. So all local functions restrict to $g_X|_C$, which proves (ii).

Since the cover is finite, there is uniform costant $v>0$ and $\zeta>0$ such that \[
\lambda_1(B_j)\geq v,\qquad \lambda_1(B_j)\lambda_2(B_j)\geq 1+\zeta. 
\]Hence \[
1-P_\chi(B_j+\chi)=\frac{\lambda_1(B_j)\lambda_2(B_j)-1}{(1+\lambda_1(B_j))(1+\lambda_2(B_j))}\geq \frac{\zeta}{1+\zeta}\left(1+\frac{1}{v}\right)^{-2}
\]This proves (iii). 

\end{proof}

\begin{remark}
    Roughly speaking, if $C$ is a smooth curve, then the above construction of $\psi$ is just a large multiple of  $\text{dist}(x,C)$, which can be used to construct cone metric near $C$.
\end{remark}
\begin{lemma}
\label{lem:finite-richberg}
For some neighborhood $U$ of $C$, there is a current
\[
 R_1\in r_X|_U
\]
with analytic singularities exactly along $C$, such that
\[
 R_1>0,\qquad R_1^2-\chi^2>0
 \quad\text{on }U\setminus C.
\]
The divisorial part of $\mu^*R_1$ on $Y$ is exactly $\tau H$.
\end{lemma}

\begin{proof}
Write $G=\theta_X+\ddc\varphi_G$ using a fixed smooth representative $\theta_X\in r_X$ and one global
quasi-psh function $\varphi_G$.  We now carry out the Richberg
construction.

\smallskip
\noindent\emph{1. Nested sets and boundary cutoffs.}
Choose $ W'_j\Subset U_j\Subset W_j$.  Put
\[
 M_j:=C\cap\partial U_j.
\]
Choose a smooth function $0\leq\vartheta_j\leq1$ on $W_j$ which
equals $1$ on a neighborhood of $M_j$ and equals $0$ on a
neighborhood of $\overline{W'_j}$.  This is possible because
$\overline{W'_j}\Subset U_j$, whereas $M_j\subset\partial U_j$.
Since there are finitely many cutoffs, there is $L>0$ such that
\begin{equation}\label{eq:cutoff-hessian-bound}
 -L\chi\leq\ddc\vartheta_j\leq L\chi
 \qquad\text{for every }j.
\end{equation}

Set
\[
 \widetilde\psi_j:=\psi_j-\eta_R\vartheta_j,\qquad
 \widetilde B_j:=G+\ddc\widetilde\psi_j,
\]
where $\eta_R>0$ will now be chosen.  This is the only place where a
cutoff Hessian enters:
\begin{equation}\label{eq:richberg-cutoff-error}
 \widetilde B_j=B_j-\eta_R\ddc\vartheta_j.
\end{equation}
Let $a,\kappa$ be as in
\eqref{eq:local-branch-uniform-margins}, decreasing $\kappa$ if
necessary so that $0<2\kappa<1$.  Choose $\eta_R$ so small that
\begin{equation}\label{eq:richberg-eta-choice}
 \eta_RL\leq\frac a2,\qquad
 \frac{\eta_RL}{a+1}\leq\frac{\kappa}{2}.
\end{equation}
The first inequality gives
\[
 \widetilde B_j\geq\frac a2\chi.
\]
For the cone estimate, put $A_j=B_j+\chi$.  Since
$A_j\geq(a+1)\chi$, equations
\eqref{eq:cutoff-hessian-bound} and
\eqref{eq:richberg-eta-choice} give
\[
 \widetilde B_j+\chi
 \geq
 \left(1-\frac{\eta_RL}{a+1}\right)A_j.
\]
Hence
\[
\begin{aligned}
 P_\chi(\widetilde B_j+\chi)
 &\leq
 \left(1-\frac{\eta_RL}{a+1}\right)^{-1}
 P_\chi(A_j)\\
 &\leq
 \left(1-\frac{\kappa}{2}\right)^{-1}(1-2\kappa)
 <1-\kappa.
\end{aligned}
\]
Thus every lowered branch is still positive and remains a strict
square-cone branch.  

\smallskip
\noindent\emph{2. Why a branch becomes irrelevant at its boundary.}
Fix $x\in M_j$.  Since the $W'_k$ cover $C$, there is an index
$k$ with $x\in W'_k$.  Near $x$,
\[
 \vartheta_j=1,\qquad\vartheta_k=0,
\]
and Lemma~\ref{lem:compatible-local-branches}(ii) gives
\[
 \widetilde\psi_j(x)-\widetilde\psi_k(x)=-\eta_R.
\]
By continuity, on a neighborhood of $x$ in $X$,
\begin{equation}\label{eq:boundary-branch-gap}
 \widetilde\psi_j
 \leq\widetilde\psi_k-\frac{\eta_R}{2}.
\end{equation}
The compact set $M_j$ is covered by finitely many such
neighborhoods.  Hence \eqref{eq:boundary-branch-gap}, with a possibly
different $k$ from one neighborhood to another, holds on a whole
collar of $M_j$.  
After shrinking one ambient neighborhood $\mathcal V$ of $C$, we
may and do assume that, for every $j$,
\[
 \mathcal V\cap\partial U_j
\]
is contained in those collars.  We also take
$\mathcal V\subset\bigcup_jW'_j$, so at least one index remains
active at every point of $\mathcal V$.

\smallskip
\noindent\emph{3. The finite regularized maximum.}
Choose $0<\delta_R<\eta_R/8$.  Let
$\rho_{\delta_R}$ be a nonnegative even radial mollifier,
supported in $(-\delta_R,\delta_R)$ and of total mass one.  For every
nonempty finite index set $I$, define
\begin{equation}\label{eq:finite-regularized-maximum}
 \mathcal M_{\delta_R,I}(t)
 :=
 \int_{\mathbb R^I}
 \max_{j\in I}(t_j+s_j)
 \prod_{j\in I}\rho_{\delta_R}(s_j)\,ds_j.
\end{equation}
This function is smooth, convex, nondecreasing in every variable, and
translation invariant:
\[
 \mathcal M_{\delta_R,I}(t_1+b,\ldots,t_{|I|}+b)
 =
 \mathcal M_{\delta_R,I}(t)+b.
\]
If one entry $t_j$ lies more than $2\delta_R$ below another entry,
then $t_j+s_j$ can never realize the maximum in
\eqref{eq:finite-regularized-maximum}.  The expression is therefore
independent of $t_j$.  Because the mollifier is a product, integrating
out $s_j$ gives exactly the same formula with the index $j$
deleted.

For $x\in U:=\mathcal V$, let
\[
 I(x):=\{j:x\in U_j\}
\]
and set locally
\begin{equation}\label{eq:richberg-global-correction}
 \psi_R(x)
 :=
 \mathcal M_{\delta_R,I(x)}
 \bigl((\widetilde\psi_j(x))_{j\in I(x)}\bigr).
\end{equation}
When $x$ crosses $\partial U_j$, inequality
\eqref{eq:boundary-branch-gap} and
$2\delta_R<\eta_R/2$ show that the $j$-th entry is irrelevant.
Thus the formulas on the two sides of $\partial U_j$ agree on a
collar.  Consequently $\psi_R$ is one globally defined smooth
function on a neighborhood of $C$.

For completeness, differentiate
\eqref{eq:richberg-global-correction}.  Put
\[
 \lambda_j:=\frac{\partial\mathcal M_{\delta_R,I}}
                   {\partial t_j}.
\]
Monotonicity and translation invariance give
\[
 \lambda_j\geq0,\qquad\sum_{j\in I}\lambda_j=1.
\]
Convexity says that the matrix
$(\mathcal M_{jk})$ of second derivatives is
positive semidefinite.  Hence
\begin{align}
 R_1
 &:=
 \theta_X+\ddc(\varphi_G+\psi_R) \notag\\
 &=
 \sum_{j\in I}\lambda_j\widetilde B_j
 +
 \sum_{j,k\in I}\mathcal M_{jk}\,
 d\widetilde\psi_j\wedge d^c\widetilde\psi_k.
 \label{eq:richberg-formula}
\end{align}
The last term is semipositive.  Thus $R_1\geq(a/2)\chi$ on $U\setminus C$.  Moreover,
\[
 R_1+\chi
 =
 \sum_j\lambda_j(\widetilde B_j+\chi)+Q,
 \qquad Q\geq0.
\]
Convexity monotonicity of $P_\chi(\cdot)$ gives
\[
 P_\chi(R_1+\chi)
 \leq
 \sum_j\lambda_jP_\chi(\widetilde B_j+\chi)
 <1-\kappa \quad \text{ on }U \setminus C.
\]
Equivalently, $R_1^2-\chi^2>0$ on $U\setminus C$.

Finally,
\[
 \varphi_G+\psi_R
\]
differs from $\varphi_G$ by a smooth function.  Hence it has exactly
the same analytic singularities as $G$, and adding this smooth
function changes no divisorial coefficient on the fixed resolution.
Lemma~\ref{lem:smaller-model} therefore gives the divisorial part
$\tau H$. 
\end{proof}

\begin{proof}[Proof of proposition~\ref{prop:analytic-descent}]
    If $C=\varnothing$, take $R=R_0=\widehat R$; the conclusion follows
from Lemma~\ref{lem:raw-pushdown}.  Assume $C\neq\varnothing$.
Let $R_1$ be the current of Lemma~\ref{lem:finite-richberg}, defined
on a neighborhood $U$ of $C$.

\smallskip
\noindent\emph{1. Strict comparison of the singular coefficients.}
Let $P_1,\ldots,P_N$ be the prime components of $F$, and write
\[
 H=\sum_{j=1}^Nh_jP_j,\qquad
 Q=\sum_{j=1}^Ne_jP_j.
\]
Here $h_j,e_j$ are nonnegative integers and
$qh_j+pe_j>0$, because $F=qH+pQ$ and
$\operatorname{Supp}F=\operatorname{Supp}H$.  Let
$\varphi_0,\varphi_1$ be potentials of $R_0,R_1$ relative to the
same smooth representative of $r_X$.  If $s_j$ is a local defining
function for $P_j$, Lemmas~\ref{lem:raw-pullback} and
\ref{lem:finite-richberg} give
\begin{align}
 \mu^*(\varphi_1-\varphi_0)
 &=
 \sum_{j=1}^N
 \bigl(\tau h_j-s(qh_j+pe_j)\bigr)\log|s_j|^2+O(1)
 \notag\\
 &=
 -s\sum_{j=1}^N
 \left(\frac q2h_j+pe_j\right)\log|s_j|^2+O(1).
 \label{eq:strict-divisorial-comparison}
\end{align}
The $O(1)$ term is locally bounded.  Every coefficient in the last
sum is positive.  Since $\operatorname{Supp}F=\mu^{-1}(C)$,
properness of $\mu$ implies
\begin{equation}\label{eq:potential-dominance-near-C}
 \varphi_1-\varphi_0\longrightarrow+\infty
 \qquad\text{as }x\longrightarrow C.
\end{equation}
Indeed, otherwise a sequence $x_\ell\to C$ with bounded difference
could be lifted to $y_\ell\in Y$; a convergent subsequence would tend
to a point of $\mu^{-1}(C)=\operatorname{Supp}F$, contradicting
\eqref{eq:strict-divisorial-comparison}.

\smallskip
\noindent\emph{2. The final regularized maximum.}
Fix $\delta_G>0$, we choose a neighborhood $V\Subset U$ of $C$. Since $\varphi_0,\varphi_1$ are smooth on $\overline{U}\setminus \overline{V} $, subtracting a sufficient large constant $M>0$, we have \[
\varphi_0>\varphi_1-M+3\delta_G=:\varphi_1'+3\delta_G, \qquad \text{ on } \overline{U}\setminus \overline{V} .
\] 
By \eqref{eq:potential-dominance-near-C}, there is another neighborhood $V_0\Subset V$ such that  \[
\varphi_1'>\varphi_0+3\delta_G \qquad \text{ on } V_0\setminus C
\]

Hence the global potential as the regularized maximum $
 \varphi:=M_{\delta_G}(\varphi_0,\varphi_1')$ is a smooth function over $X\setminus C$, which equals $\varphi_0$ outside $\overline{V}$ and equals $\varphi_1'$ near $C$. Set $R=\theta_X+\ddc \varphi$. From the convexisty of $P_\chi$, $R$ satisfies the strict cone inequality over $X\setminus C$. 

\smallskip
\noindent\emph{3. The weak cone inequality across $C$.}
The analytic curve $C$ has Lebesgue measure zero.  The current just
constructed is positive, smooth off $C$, and satisfies
$R^2-\chi^2>0$ there.  Lemma~\ref{lem:convolution-negligible-set}
therefore gives
\[
 R_\delta^2-\chi_0^2\geq0
\]
for every admissible constant background $\chi_0$. 
\end{proof}

\begin{proof}[Proof of Theorem~\ref{thm:main-thm}]
Lemma~\ref{lem:rho-nef} gives finitely many $J$-null curves, and
lemma~\ref{lem:finite-null-divisor} gives  finitely many $J$-null prime divisors.
Their union is precisely $\JNull(\alpha,\beta)$, proving (i).

Let $R$ be the current
of proposition~\ref{prop:analytic-descent}, and put
\[
 A:=\frac{R+\chi}{\cab},\qquad
 \sigma:=\frac{A_0s^2}{\cab},\qquad
 S:=A+\sigma[E],\qquad
 T:=S+\varepsilon\chi.
\]
Then
\[
 [S]=\alpha-\varepsilon\beta,\qquad S\geq0,
\]
and $P_\chi(S)\leq \cab$ in the local-convolution sense:   On a coordinate ball write $R=\ddc u$ and
$\chi=\ddc q$.  Let $\chi_0$ be any constant-coefficient K\"ahler
form with $\chi\geq\chi_0$, and use a subscript $\delta$ for
ordinary convolution on a smaller ball.  Proposition~\ref{prop:analytic-descent}
gives
\[
 R_\delta^2-\chi_0^2\geq0.
\]
Set
\[
 A_\delta^0:=\frac{R_\delta+\chi_0}{\cab}.
\]
An ordinary smooth calculation gives
\[
 \cab(A_\delta^0)^2-2\chi_0\wedge A_\delta^0
 =\frac{1}{\cab}(R_\delta^2-\chi_0^2)\geq0.
\]
On the other hand, the convolution of $A=(R+\chi)/c$ is
\[
 A_\delta=A_\delta^0+\frac{\chi_\delta-\chi_0}{\cab}
 \geq A_\delta^0,
\]
because convolution preserves the inequality $\chi\geq\chi_0$.
So $A_\delta$ still
satisfies the same inequality.  Finally, convolving the positive
divisor current $\sigma[E]$ only adds another semipositive matrix.
Thus the convolution of $S=A+\sigma[E]$ satisfies
\[
 cS_\delta^2-2\chi_0\wedge S_\delta\geq0,
\]
which is precisely $P_\chi(S)\leq c$.

It remains to show the analytic singularities. Let $\cI=(f_1,...,f_N)$ br the coherent ideal used in construction of $R$ near $C$, i.e. $\cI\cdot \cO_Y=\cO(-H)$.  Let $g_i$ be a local defining function for $D_i$. A local potential of $S$ has singularity part \[
\frac{1}{\cab}\left[\tau\log\left(\sum_{\ell}|f_\ell|^2\right)+A_0s^2\sum_ia_i\log|g_i|^2 \right].
\]
Since $\tau, A_0, s, a_i$ are all rational, we can choose $N_0\in \N$ such that \[
m_0:=N_0\tau\in \N,\qquad m_i:=N_0A_0s^2a_i\in \N
\]for every $i$. Then singular part is now \[
\frac{1}{\cab N_0}\log\left[\left(\sum_\ell|f_\ell|^2\right)^{m_0}\prod_{i}|g_i|^{2m_i}\right].
\]
One can write the expression in the bracket as a finite sum of squared absolute values of the generators of the coherent product ideal \[
\cI^{m_0}\cdot\cO_X\left(-\sum_im_iD_i\right).
\]Thus $S$ has analytic singularity as desired. 

Finally, proposition~\ref{prop:easy-direction} gives $\Null(\alpha,\beta)\subset E_+(T')$ for every $T'\in \cK_{J}^{\rm an}(\alpha,\beta)$. Hence $\Null(\alpha,\beta)\subset \JEnK(\alpha,\beta)\subset E_+(T)=\Null(\alpha,\beta)$. This proves the asserted equality. 
\end{proof}

\begin{rmk}[What has actually been regularized]
The proof does not assert that Demailly regularization preserves
$R^2-\chi^2\geq0$.  The raw pushdown $R_0$ already has the weak cone
property by convexity of the constant-background matrix cone.  Its
singularities are replaced near $C$ by a neat model with the smaller
coefficient $\tau H$, while the raw coefficient is
$sF=s(qH+pQ)$.  The strict componentwise comparison
$\tau h_j<s(qh_j+pe_j)$ follows from $\tau=qs/2$; the replacement
is made by a regularized maximum only where both inputs are smooth and
strict.
\end{rmk}

\section{Convergence of the $J$-flow outside the J-null locus}
Let $Z=\JNull(\alpha,\beta).$
The proof of the main theorem gives a current \[
T=S+\epsilon\chi\in \alpha, \qquad, S\ge 0, \qquad [S]=\alpha-\varepsilon\beta,
\]
with analytic singularitiese exactly along $Z$ and with $P_\chi(S)\leq \cab$ in the local-convolution sense. It follows that there are constants $b_0,\sigma>0$ such that \[
T\geq b_0\chi,\qquad P_\chi(T)\le c-\sigma.
\]
We fix a K\"ahler metric $\omega_0\in \alpha$. There is $\rho\in \Psh(X,\omega_0)$ with analytic singularity along $Z$ such that $T=\omega_0+dd^c\rho$. We normalize $\rho$ such that $\rho\leq 0$. In particular, we have $\rho\to -\infty$ along $Z$. 

Recall the $J$-flow \eqref{eq:J-flow},
\begin{equation}
 \dot\varphi=\cab-\tr_{\omega_\varphi}\chi,\qquad
 \omega_\varphi=\omega_0+\ddc\varphi>0.
\end{equation}
The smooth solution exists for long time
\cite[Theorem~1.1]{C04}.

Put $v=\dot\varphi$.  Differentiating \eqref{eq:J-flow} gives
\begin{equation}\label{eq:v-equation}
 v_t=h^{i\bar j}v_{i\bar j},\qquad
 h^{i\bar j}=(\omega_\varphi)^{i\bar p}
       (\omega_\varphi)^{q\bar j}\chi_{q\bar p}.
\end{equation}
The maximum principle gives
\begin{equation}\label{eq:v-bound}
 \|v(t)\|_{L^\infty(X)}
 \leq\|v(0)\|_{L^\infty(X)}.
\end{equation}
Since $\tr_{\omega_\varphi}\chi=\cab-v$, every eigenvalue of
$\chi$ relative to $\omega_\varphi$ is bounded by their bounded
sum.  Hence
\begin{equation}\label{eq:metric-lower}
 \omega_\varphi\geq\kappa\omega_0
\end{equation}
for a uniform $\kappa>0$.
\begin{lemma}\label{lem:v-decay}
    Along the flow \eqref{eq:J-flow}, we have $\dot\varphi\to 0$ 
in $L^2(\omega_0^3)$. \end{lemma}
\begin{proof}
For $\tau>0$, set $\chi_\tau=\chi+\tau\omega_0$, and $\delta_\tau=\tau/2$, $c_\tau=\cab+3\tau$. Then for every irreducible $p$-dimensional subvariety, $p=1,2$, \[
\int_V((c_\tau-(3-p)\delta_\tau)\alpha^p-p\chi_\tau\alpha^{p-1}=\int_V(\alpha^p-p\chi\alpha^{p-1})+(3-p)(\tau-\delta_\tau)\int_V\alpha^p\ge0.
\]
Therefore, by Chen's theorem\cite[Theorem 1.1]{C21}, It is  equivalent to coercivity of the $\cJ_{\chi_\tau}$ functional. In particular, $\cJ_{\chi_\tau}$ is bounded below. On the other hand, along the flow \eqref{eq:J-flow}, 
 \begin{equation}
 \begin{aligned}
     \frac{d}{dt}\left(\cJ_\chi(\varphi)-\cJ_{\chi_\tau}(\varphi)\right)=&\int_X\dot\varphi(3\chi\omega_\varphi^2-\cab\omega_\varphi^3-3\chi_\tau\omega_\varphi^2+c_\tau\omega_\varphi^3)\\
     =&\int_X \dot\varphi(3\tau\omega_0\omega_\varphi^2-3\tau\omega_\varphi^3)
 \end{aligned}
\end{equation}
By \eqref{eq:v-bound}, we have \[
 \left|\frac{d}{dt}\left(\cJ_\chi(\varphi)-\cJ_{\chi_\tau}(\varphi)\right)\right|\le C\tau.
\]
Thus \begin{equation}
    \label{eq:pert-j-function-estimate}\cJ_{\chi}(\varphi)\geq-C\tau t-C_\tau.
\end{equation}
The $\cJ_\chi$ functional is convex along the $J$-flow: \[\begin{aligned}
  \frac{d^2}{dt^2}\cJ_{\chi}(\varphi)=&\frac{d}{dt}\int_X-v^2\omega_\varphi^3\\
  =&-2\int_Xv\dot v\omega_\varphi^3-\int_Xv^2\Delta_{\varphi}v\omega_\varphi^3\\
  =&-2\int_X vg_\varphi^{i\bar p}g_\varphi^{q\bar j}\chi_{q\bar p}v_{i\bar j}\omega_\varphi^3-+2\int_X|\nabla v|_{\omega_\varphi}^2v\omega_\varphi^3\\
=&2\int_Xg_\varphi^{i\bar p}g_\varphi^{q\bar j}\chi_{q\bar p} v_iv_{\bar j}\omega_\varphi^3\geq 0,
\end{aligned}\]
where we have used in the last line that \[
g_\varphi^{i\bar p}g_\varphi^{q\bar j}\chi_{q\bar p,i}v_{\bar j}=g_\varphi^{i\bar p}g_\varphi^{q\bar j}\chi_{i\bar p,q}v_{\bar j}=g_\varphi^{q\bar j}(g_{\varphi}^{i\bar p}\chi_{i\bar p})_{q}v_{\bar j}=|\nabla v|_{\omega_\varphi}^2.
\]
Hence $\frac{d}{dt}\cJ_\chi(\varphi)$ is nondecreasing, and the $\lim_{t\to \infty}\frac{d}{dt}\cJ_\chi(\varphi)$ exists, let us denote it by $-L\leq 0$.
It follows that \[
\lim_{t\to \infty}\frac{\cJ_{\chi(\varphi)}}{t}=\lim_{t\to \infty}\frac{d}{dt}\cJ_\chi(\varphi)=-L.
\]
From \eqref{eq:pert-j-function-estimate}, we get $-L\geq -C\tau.$ Since $\tau>0$ is arbitrary, we get $L=0$. 

\eqref{eq:metric-lower} and $\frac{d}{dt}\cJ_{\chi}(\varphi)=-\int_{X}v^2\omega_\varphi^3$ give $v\to 0$ in $L^2(\omega_0^3)$.
\end{proof}

Now we prove the relative $C^0$-estimate. 
We use the following standard parabolic ABP lemma.

\begin{lemma}\label{lem:ABP}
Let $U\subset B_r(0)\subset\mathbb R^m$ be a bounded
domain containing $0$, and let $Q=U\times[-1,0]$.
Write
\[
 \partial_PQ
 =\bigl(\partial U\times[-1,0]\bigr)
  \cup\bigl(\overline{U}\times\{-1\}\bigr).
\]
 Suppose
 $W\in C^\infty(\overline Q)$ and
\[
 W(0,0)=L,\qquad W\geq L+h\quad\text{on }\partial_PQ.
\]
There is a measurable set $\Gamma\subset Q$ on which
\[
 D_x^2W\geq0,\qquad W_t\leq0,\qquad
 W\leq L+\frac{5h}{8},
\]
and
\begin{equation}\label{eq:ABP-general}
 c_m\frac{h^{m+1}}{r^m}
 \leq
 \int_\Gamma(-W_t)\det_{\mathbb R}(D_x^2W)\,dx\,dt.
\end{equation}
\end{lemma}

\begin{lemma}\label{lem:C0-estimate}
    There is a uniform constant $C$ such that \[
    \varphi-\bar{\varphi}-\rho\geq -C \quad \text{ on } X\setminus Z\times [0,\infty) \quad \text{where }\bar\varphi=\frac{\int_{X}\varphi(t)\chi^3}{\int_X\chi^3}.
    \]
\end{lemma}

\begin{proof}
    Set $w=\varphi-\bar{\varphi}-\rho$, $B=T|_{X\setminus Z}$. Then $\omega_\varphi=B+dd^cw$. 
On every finite time interval, $w\to+\infty$ uniformly as approaching 
$Z$.

Choose $t_0$ so large that
\[
\big|\frac{d}{dt}\bar{\varphi}\big|\leq\sigma/16\qquad(t\geq t_0),
\]which is possible since $v\to 0$ in $L^2(\omega_0^3)$. 

The estimate on a fixed finite interval is automatic.  Suppose that,
for some $T>t_0+2$, the minimum of $w$ on
$U\times[t_0,T]$ is $L\ll-1$, attained at
$(x_*,t_*)$ with $t_*\geq t_0+2$.  Work in a fixed-radius
coordinate ball $B_{2r}$ for $\chi$, centered
at $x_*=0$, on which $dd^c(|z|^2/r^2)\leq C_r\chi$. Now we choose constant $\eta>0$ such that \begin{equation}\label{eq:eta-choice}
 2\eta\leq\frac{3\sigma}{16},\qquad
 C_r\eta\leq
 \min\left\{\frac{b_0}{2},\frac{\sigma b_0^2}{32}\right\}.
\end{equation}
Set \[
I_*=[t_*-1.t_*], \qquad q(z,t)=\eta\left(\frac{|z|^2}{r^2}+(t-t_*)^2\right),
 \qquad W=w+q.
\]
Because $L$ is the global spacetime minimum,
\[
 W(0,t_*)=L,\qquad W\geq L
 \quad\text{on }(B_r\setminus Z)\times I_*.
\]
Let $M>0$ be a sufficiently large constant, and let $U_M:=\text{ the connected component containing $x_*$ of }\{\rho>-M\}\cap B_{r}$, and $m_*:=\min_{X\times I_*}(\varphi-\bar\varphi).$ Then \[
\partial U_M\subset (\partial B_r\cap\{\rho\geq -M\})\cup(B_r\cap\{\rho=-M\}).
\]
We now check every part of the parabolic boundary explicitly.  If
$z\in\partial B_r\cap\partial U_M$, then
\[
 W(z,t)\geq L+q(z,t)\geq L+\eta
 \qquad(t\in I_*).
\]
If $z\in B_r\cap\partial U_M$, then
$\rho(z)=-M$, and hence
\[
 W(z,t)=u(z,t)+M+q(z,t)
 \geq m_*+M\geq L+\frac\eta2
 \qquad(t\in I_*).
\]
Finally, on the bottom face $t=t_*-1$,
\[
 W(z,t_*-1)\geq L+\eta
 \qquad(z\in\overline{U_M}).
\]
Define $\widetilde W(z,s)=W(z,t_*+s)$ for $-1\leq s\leq0$.
Define $\widetilde W(z,s)=W(z,t_*+s)$ for $-1\leq s\leq0$.
The preceding three inequalities say exactly that all the hypotheses
of Lemma~\ref{lem:ABP} hold for $\widetilde W$ on
$U_M\times[-1,0]$, with $m=6$ and $h=\eta/2$.  The
constant in that lemma depends only on $r$, independent of $M$.  We obtain
\begin{equation}\label{eq:ABP-applied}
 C_0\eta^7
 \leq
 \int_{\Gamma_W}(-W_t)\det_{\mathbb R}(D_x^2W)\,dx\,dt.
\end{equation}
On $\Gamma_W$,
\begin{equation}\label{eq:contact-data}
 D_x^2W\geq0,\qquad W_t\leq0,\qquad
 W\leq L+\frac{5\eta}{16},
\end{equation}
and
\begin{equation}\label{eq:contact-complex}
 \ddc w\geq-C_r\eta\chi,\qquad
 \varphi_t-(\bar{\varphi})_t=w_t\leq2\eta.
\end{equation}
We claim that both $B$ and $\omega_\varphi$ are uniformly bounded at every point of $\Gamma$, although $B$ is not bounded a priori.    Choose normal coordinates of $\chi$, diagonalizing $\omega_\varphi$, with eigenvalues $\lambda_i$. Since $\omega_\varphi=B+dd^c w$ and \eqref{eq:contact-complex}, $B\leq \omega_\varphi+C_r\eta\chi$. The flow equation gives \begin{equation}\label{eq:sum-eigenvalues}
\sum_{i=1}^3\frac{1}{\lambda_i}=\cab-w_t-(\bar{\varphi})_t\geq \cab-2\eta-(\bar{\varphi})_t\geq\cab-\frac{\sigma}{4},
\end{equation}where we used $|(\bar{\varphi})_t|\leq \sigma/16$, $w_t\leq 2\eta$, and \eqref{eq:eta-choice}. 
Since $B\geq b_0\chi$, and $B\leq \omega_\varphi+C_r\eta\chi$, we have $\lambda_i\geq b_0-C_r\eta\geq \frac{b_0}{2}, $ and $P_\chi(\omega_\varphi+C_r\eta\chi)\leq P_\chi(B)\leq\cab-\sigma.$ 

Since $\lambda_j\geq b_0/2$, \[
0\leq \frac{1}{\lambda_j}-\frac{1}{\lambda_j+C_r\eta}
=\frac{C_r\eta}{\lambda_j(\lambda_j+C_r\eta)}\leq \frac{4C_r\eta}{b_0^2}.\]
For $C_r\eta\leq \frac{\sigma b_0^2}{32}$, it follows that \[
\sum_{i\neq j}\frac{1}{\lambda_j}\leq \sum_{i\neq j}\left(\frac{1}{\lambda_j+C_r\eta}+\frac{\sigma}{8}\right)\leq P_\chi(\omega_\varphi+C_r\eta\chi)+\frac{\sigma}{4}\leq \cab-\frac{3\sigma}{4}.
\]
Combined with \eqref{eq:sum-eigenvalues}, we get $\frac{1}{\lambda_i}\geq \frac{\sigma}{2}$ for every $i$. The upper bound for $\omega_\varphi$ follows which confirms the claim. 

It follows from $dd^cW=\omega_\varphi-B+dd^cq$ that $\tr_\chi(dd^c W)\leq C$ for a unoform constant $C$ depending only on $b_0,\sigma,r$. It gives $\Delta_{\R}W=2\sum_iW_{i\bar i}\leq C$. Sicne $D^2_xW\geq 0$, all its real eigenvalues are nonnegative and therefore bounded by this real trace. Hence $\det_{\mathbb R}(D_x^2W)\leq C$. Also,  \[\varphi_t-(\bar{\varphi})_t=\cab-\tr_{\omega_\varphi}\chi\geq -C, \qquad q_t\geq -2\eta.\]
Together with $W_t\leq 0$, this yields $0\leq -W_t\leq C$. Consequently \eqref{eq:ABP-applied} implies
\begin{equation}\label{eq:contact-positive-volume}
 \operatorname{Vol}(\Gamma_W)\geq C_1>0.
\end{equation}

On the other hand, $\rho\leq0$ gives $w=\varphi-\bar{\varphi}-\rho\geq \varphi-\bar{\varphi}$, and hence
$w_-\leq (\varphi-\bar{\varphi})_-$.  If $|L|\geq2\eta$, then
\eqref{eq:contact-data} gives $w^-\geq|L|/2$ on $\Gamma_W$.

 Since $\varphi-\bar\varphi\in \Psh(X,\omega_0)$ with zerom mean, by Green's identity, there is a uniform constant $C$ such that $\int_X|\varphi-\bar\varphi|\chi^3\leq C. $
 It follows that 
\[
 \operatorname{Vol}(\Gamma_W)
 \leq\frac{2}{|L|}\int_{\Gamma_W}w^-
 \leq\frac{C}{|L|}.
\]
This contradicts \eqref{eq:contact-positive-volume} for large
$|L|$.  The bound is independent of $T$, so $T\to\infty$
proves the lemma.

\end{proof}

Now we prove higher order estimates off $Z$. 

\begin{lemma}\label{lem:local-c2}
    For every $K\Subset U$ and $k\geq0$, there are constants
$C_K,C_{K,k}$ such that, for all sufficiently large $t$,
\[
 C_K^{-1}\chi\leq\omega_\varphi(t)\leq C_K\chi,
 \qquad
 \|\varphi-\bar\varphi\|_{C^k(K)}\leq C_{K,k}.
\]
\end{lemma}

\begin{proof}
    Condsider $Q=\log\tr_\chi\omega_\varphi-Aw$, where $w=\varphi-\bar{\varphi}-\rho$. At any time slice, $Q(x,t)\to \infty$ as $x\to Z$. Write $\Lambda=\tr_\chi\omega_\varphi$. The Song-Weinkove calculation \cite[Lemma 3.1]{SW08}, gives \[
    (h^{i\bar j}\partial_i\partial_{\bar{j}}-\partial_t)\log \Lambda\geq \frac{1}{\Lambda}\left(
 h^{k\bar\ell}R_{k\bar\ell}{}^{i\bar j}
                  (\omega_\varphi)_{i\bar j}
 -(\omega_\varphi)^{k\bar\ell}R_{k\bar\ell}
 \right),
    \]
    Where $R$ is the curvarure of $\chi$. The right-hand side is bounded below by $-C$, using only
$\omega_\varphi\geq\kappa\chi$: in $\chi$-unitary
coordinates, $|h|_\chi\leq\kappa^{-2}$,
$|\omega_\varphi|_\chi=\Lambda$,
$|\omega_\varphi^{-1}|_\chi\leq\kappa^{-1}$, and
$\Lambda\geq3\kappa$. Furthermore, $\ddc w=\omega_\varphi-B$ gives
\begin{equation}\label{eq:SW-curvature}
 (h^{i\bar j}\partial_i\partial_{\bar j}-\partial_t)w
 =
 2\tr_{\omega_\varphi}\chi
 -h^{i\bar j}B_{i\bar j}-\cab+(\bar{\varphi})_t.
\end{equation}
At a positive-time maximum of $Q$, the parabolic maximum principle,
\eqref{eq:SW-curvature}, and the last identity give
\[
 0\geq(h^{i\bar j}\partial_i\partial_{\bar j}-\partial_t)Q
 \geq
 A\bigl(\cab-(\bar\varphi)_t+h^{i\bar j}B_{i\bar j}
             -2\tr_{\omega_\varphi}\theta\bigr)-C.
\]
Thus, choosing $A$ sufficiently large,
\begin{equation}\label{eq:max-evolution}
 \cab-(\bar\varphi)_t+h^{i\bar j}B_{i\bar j}
 -2\tr_{\omega_\varphi}\theta
 \leq\varepsilon_0,
\end{equation}
where $\varepsilon_0>0$ is as small as desired.
At that point, take normal coordinates of $\chi$ in which
$\omega_\varphi=\operatorname{diag}(\lambda_i)$, and write
$b_i=B_{i\bar i}$.  Then
\[
 h^{i\bar j}B_{i\bar j}=\sum_i\frac{b_i}{\lambda_i^2},
 \qquad
 \tr_{\omega_\varphi}\chi=\sum_i\frac1{\lambda_i}.
\]
Increase $t_0$ so that $|(\bar\varphi)_t|\leq\sigma/4$, and choose
$\varepsilon_0\leq\sigma/4$. 

Observe that $\sup_j\sum_{i\neq j}\frac{1}{b_i}\leq P_\chi(B)\leq \cab-\sigma$: In deed, let $\nu_1\ge\nu_2\ge\nu_3>0$ be the eigenvalues of $B^{-1}$ and $e_i$ be the $i$-th vector of normal coordinates of $\chi$. Then $b_i=\langle Be_i,e_i\rangle_{\chi}$ and $(B^{-1})_{j\bar j}\geq \nu_3$ by Rayleigh quotient. So \[
1=|\langle B^{1/2}e_i,B^{-1/2}e_i\rangle_\chi|^2\leq \|B^{1/2}e_i\|_\chi^2\|B^{-1/2}e_i\|_\chi^2=\langle Be_i,e_i\rangle_\chi\langle B^{-1}e_i,e_i\rangle_\chi=b_i(B^{-1})_{i\bar i}.
\]
Consequently, for fixed $j$,
\[
\begin{aligned}
\sum_{i\ne j}\frac1{b_i}
&\le \sum_{i\ne j}(B^{-1})_{i\bar i}\\
&=\operatorname{tr}(B^{-1})-(B^{-1})_{j\bar j}\\
&\le \nu_1+\nu_2+\nu_3-\nu_3\\
&=P_\chi(B)\le\cab-\sigma,
\end{aligned}
\]

Now, for $i\ne j$, \[
\frac{b_i}{\lambda_i^2}-\frac{2}{\lambda_i}=b_i\left(\frac{1}{\lambda_i}-\frac{1}{b_i}\right)^2-\frac{1}{b_i}.
\]
It follows that for any $j$, \[
\cab-(\bar{\varphi})_t+\sum_{i}\frac{b_i}{\lambda_i^2}-2\sum_{i}\frac{1}{\lambda_i}\geq \frac{3\sigma}{4}-\frac{2}{\lambda_j}
\]
Combined with \eqref{eq:max-evolution}, we get $\lambda_j\leq \frac{4}{\sigma}$. Hence $\tr_\chi\omega_\varphi\leq \frac{12}{\sigma}$ at the maximum of $Q$. 
If the maximum lies on $t=t_0$, the same
conclusion follows from the fixed initial-time metric.  Comparing
$Q$ with its maximum yields
\begin{equation}\label{eq:weighted-trace}
 \tr_\chi\omega_\varphi
 \leq
 C\exp\!\left(
 A\left[w-\inf_{U\times[t_0,t]}w\right]\right).
\end{equation}
Lemma~\ref{lem:C0-estimate} bounds the infimum.  On
$K\Subset U$, $\rho$ is bounded and $\sup_Xu\leq C$, so
\eqref{eq:weighted-trace} gives
\[
 \tr_\chi\omega_\varphi\leq C_K.
\]
Together with \eqref{eq:metric-lower}, this is the local two-sided
metric estimate.

Higher order estimates follows from the Evans-Krylov estimates and Schauder bootstrapping. 
\end{proof}

\begin{lemma}\label{lem:subsequence-weak-limit}
    Every sequence
$t_j\to\infty$ has a subsequence, still denoted $t_j$, for which
\[
( \varphi-\bar\varphi)(t_j)\longrightarrow \varphi_*
 \quad\text{in }L^1(X)\cap C^\infty_{\rm loc}(X\setminus Z),
\]
and
$\Omega_*:=\omega_0+\ddc \varphi_*$ is a global K\"ahler current satisfying
\[
 \Omega_*\geq\kappa\chi,\qquad
 \tr_{\Omega_*}\chi=\cab\quad\text{on }X\setminus Z.
\]
\end{lemma}
\begin{proof}
Quasi-psh compactness first give an
$L^1$-convergent subsequence.  Lemma~\ref{lem:local-c2} and a
diagonal extraction give smooth convergence on $X\setminus Z$, to the same
distributional limit.  The global inequality
$\omega_{\varphi(t)}\geq\kappa\chi$ passes to the weak current
limit.  Finally,
\[
 \cab-\tr_{\omega_{\varphi(t_j)}}\chi=v(t_j)
 \longrightarrow\cab-\tr_{\Omega_*}\chi
 \quad\text{smoothly locally on }x\setminus Z.
\]
The left side tends to zero in $L^2(X,\chi^3)$ by
Lemma~\ref{lem:v-decay}; hence the smooth local limit vanishes.
\end{proof}

We can identify the limit as the one obtained by Murakami \cite{M26b}: there is a unique weak current $\Omega_\infty\in \alpha$ and $\omega_{\varphi(t)}\to \Omega_\infty$ as currents. The current $\Omega_\infty$ satisfies \[
\cab\langle\Omega_\infty^3\rangle=3\chi\wedge\langle\Omega_\infty^2\rangle, \qquad \cab\langle\Omega_\infty^2\rangle-2\chi\wedge\langle\Omega_\infty\rangle\geq 0, \qquad \cab\langle \Omega_\infty\rangle-\chi\geq 0. 
\]
Therefore, the subsequence convergence of $\varphi-\overline{\varphi}$ is actually full convergence. This proves corollary \ref{cor:convergence-J-flow}.

\AtNextBibliography{\small}
\begingroup
\setlength\bibitemsep{2pt}
\printbibliography

\end{document}